\documentclass{article}
\usepackage{amssymb,amsmath}
\usepackage{graphicx}
\usepackage{color}
\usepackage{accents}

\newtheorem{theorem}{Theorem}[section]
\newtheorem{lemma}[theorem]{Lemma}
\newtheorem{remark}[theorem]{Remark}
\newtheorem{corollary}[theorem]{Corollary}

\newcommand{\id}{\mathop{\rm id}}
\newcommand{\inv}{\mathop{\rm inv}}
\newcommand{\eps}{\varepsilon}
\newcommand{\bea}{\begin{eqnarray*}}
\newcommand{\eea}{\end{eqnarray*}}
\newcommand{\be}{\begin{equation}}
\newcommand{\ee}{\end{equation}}
\newcommand{\hphi}{{\hat\phi}}

\newcommand{\cla}{{\cal A}}
\newcommand{\clb}{{\cal B}}

\newcommand{\clt}{{\cal T}}
\newcommand{\cll}{{\cal L}}
\newcommand{\cls}{{\cal S}}
\newcommand{\tr}{{\mathop{\rm tr}}}
\newcommand{\Tr}{{\mathop{\rm Tr}}}
\newcommand{\qed}{\hspace*{\fill}\rule{2mm}{2mm}}
\newcommand{\RRM}{\mathbb{R}}
\newcommand{\CCM}{\mathbb{C}}
\newcommand{\HHM}{\mathbb{H}}
\newcommand{\TTM}{\mathbb{T}}
\newcommand{\supp}{\mathop{\rm supp}}

\newcommand{\diam}{\mathop{\rm diam}}
\newcommand{\Hol}{\mathop{\rm Hol}}
\renewcommand{\Re}{\mathop{\rm Re}}
\renewcommand{\Im}{\mathop{\rm Im}}
\newcommand{\Mathring}[1]{\underaccent{\circ}{#1}}
\newcommand{\etk}{\eta^{(k)}}
\usepackage{amssymb}

\newenvironment{keywords}%
   {\begin{trivlist}\item[]{\bfseries\sffamily Keywords:}\ }% oder "Keywords:"
   {\end{trivlist}}

\newenvironment{MSC}%
   {\begin{trivlist}\item[]{\bfseries\sffamily MSC:}\ }% oder  
   {\end{trivlist}}

\numberwithin{equation}{section}

\begin{document}

\title{Well-posedness for a moving boundary model of an evaporation front in a 
porous medium}
\author{Friedrich Lippoth \\ \small{Institute of Applied Mathematics, Leibniz 
University Hannover, Welfengarten 1,} \\ \small{ D-30167 Hannover, Germany,}
 \small{\texttt{lippoth@ifam.uni-hannover.de}} \\
%\\Mark A. Peletier\\ 
\\Georg Prokert\\
\small{Faculty of Mathematics and Computer
Science, TU Eindhoven} \\ \small{P.O. Box 513 5600 MB Eindhoven, The
Netherlands,} 
\small{\texttt{g.prokert@tue.nl}}
}

\date{}
\maketitle

\begin{center}
\end{center}

\begin{abstract}
We consider a two-phase elliptic-parabolic moving boundary problem modelling
an evaporation front in a porous medium. Our main result is a proof of short-time 
existence and uniqueness of strong solutions to the corresponding nonlinear 
evolution problem in an $L^p$-setting.  It relies critically on nonstandard optimal 
regularity results for a linear elliptic-parabolic system with dynamic boundary condition.
\end{abstract}


\begin{MSC}
Primary $35$R$37$, Secondary $35$M$33$, $76$T$10$ 
\end{MSC}

\section{Introduction} \label{secintro}
The classical Stefan and Hele-Shaw problems are probably the best studied
representatives of a wide class of moving boundary problems arising from a broad
variety of models in continuum mechanics, other fields of physics as well as in
the life sciences. One of the standard techniques for a
rigorous mathematical treatment of these problems consists in transforming the
problem under consideration to a fixed reference domain by a time-dependent
diffeomorphism and to apply methods from functional analysis to the resulting
evolution problems. These problems are typically strongly nonlinear, nonlocal,
and have parabolic character. In connection with this character, a natural
well-posedness condition on the parameters and/or data occurs which often
has a direct interpretation in terms of the underlying model.

The present paper starts a discussion, along these lines, of a two-phase
problem arising from a model for flow with evaporation in a porous medium, with
gravity as driving force. The two phases
represent a porous medium whose free pore space is filled either by a liquid
(water, phase ``$-$'') or by its vapor, resulting in variable humidity (phase
``$+$''). Mathematically, this leads to an elliptic governing equation for the
liquid pressure (as in Hele-Shaw problems) in the liquid phase and a parabolic
governing equation for the humidity (as in Stefan problems) in the vapor phase.
The motion of the phase boundary (which is supposed to be a sharp interface) is
governed by conservation of mass and the fact that at fixed temperature and
vapor pressure, condensation of the vapor occurs at a certain fixed humidity. 
A remarkable point here is that the water is situated above the vapor, which 
gives rise to possible instabilities. 

This problem has been investigated from a modelling point of view, with 
emphasis on (in)stability analysis of horizontal equilibria in dependence of 
the physical parameters, by Schubert and Straus \cite{schu}, Ilichev and 
Shargatov \cite{ilshar}, and Ilichev and Tsypkin \cite{iltscat}. Our aim here 
is to give a strict short-time existence and uniqueness result for the 
nonlinear moving boundary problem and to explicitly identify the well-posedness 
condition in terms of the initial data and the dimensionless parameters. This 
condition can be viewed as a generalization to both the well-posedness 
conditions for the Hele-Shaw and the Stefan problem, as formally neclecting one 
of the phases recaptures these conditions.

Stefan-type and Hele-Shaw-type problems, both in one- and two-phase settings, 
have been studied extensively from a mathematical point of view. In problems where the 
motion of the free boundary is governed by both an elliptic and a parabolic equation in 
the bulk phase, however, most work has been devoted to surface evolutions dominated by 
a single highest-order term representing the influence of curvature, as e.g. in 
the work of Escher in a tumor model \cite{es},  
the references given there, and in \cite{lp}, where a Stokes flow problem with 
osmosis is investigated. The only exception known to us is 
a result by Bazalii and Degtyarev
\cite{baz}, who show well-posedness for short time for a coupled 
elliptic-parabolic moving boundary problem 
(with boundary conditions different from  those considered here) by means of 
parabolic regularization in a H\"older 
space setting.

As in [1], a specific difficulty arises from the fact that the boundary conditions at the 
interface do not involve curvature but normal derivatives from both the elliptic and 
parabolic phase. Because of this feature the corresponding (linear, constant-
coefficient, halfspace) model problem is nonstandard, more precisely, its 
corresponding operator symbol is inhomogeneous. To derive the necessary estimates, 
we use recently established results on parabolic problems of this 
type, systematically presented by Denk and Kaip \cite{deka}. 
Based on this, the main technical effort is in carrying over the necessary
estimates to the variable coefficient case. Again, although the basic approach
of ``freezing of coefficients'' is straightforward, we cannot rely directly on
standard results here due to the coupling of an elliptic and a parabolic
phase. Moreover, as we work in an $L^p$-setting oriented at the one used by 
Solonnikov and Frolova in
\cite{frol,solfr} for the one-phase Stefan problem, one has to work with Besov 
spaces
of ``negative differentiability in space'' in the elliptic phase, and to
exploit the parabolic character of the problem by working simultaneously with
vector-valued function spaces from the same class, but with different
smoothness parameters, see Theorem \ref{mainthrmlin}. 

The present paper is organized as follows: In the remainder of Section \ref{secintro} we 
derive our moving boundary problem (in a spatially periodic setting) from the underlying 
physical model. We explicitly include the nondimensionalization and formulate 
the  well-posedness 
condition (\ref{sign}) in this setting.  Section \ref{sectrans} is devoted to 
the transformation 
of our problem to a fixed domain and contains the formulation of our main 
result, Theorem \ref{mainthrm}. 
Section \ref{seclin} discusses a sequence of linear problems,  starting from a 
half-space model problem 
and leading up to the full linearization of the problem under consideration.  
The results on this linearization 
are applied in Section \ref{secnonlin} to prove our main well-posedness theorem. 
 The appendix contains 
a number of technical results whose proofs we include for completeness and  
convenience, without claiming 
originality.

\subsection{Problem setting}

\begin{figure}[t]
\begin{center}
\includegraphics[width=0.6\textwidth]{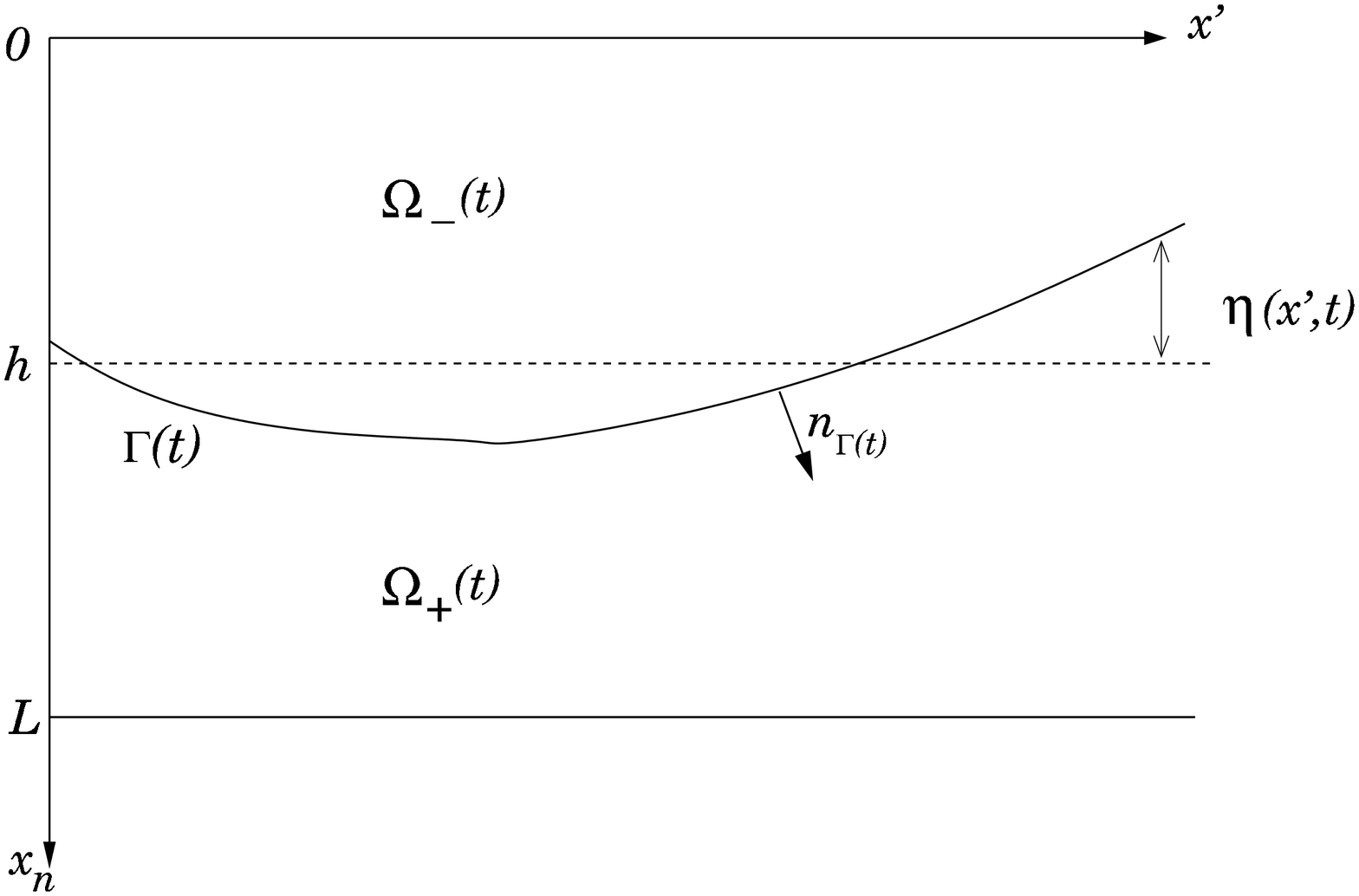}
\end{center}
\caption{\label{notation} Geometric setting. Note that the $x_n$ axis is 
oriented downwards, and that the liquid phase lies above the vapor phase.}
\end{figure}

Let $n \in \mathbb{N}$, $n \geq 2$, and let $\TTM^{n-1} := \RRM^{n-1} / \mathbb{Z}^{n-1} $ 
be the $(n-1)$-dimensional torus. We assume that the porous medium occupies a layer domain
$\Omega:=\TTM^{n-1}\times(0,L)$, with the $n$-th unit vector oriented
``downwards'', i.e. in the direction of gravity. The domain is separated in
two phases by an interface depending on time $t\in[0,T]$:
\bea
\Omega&=&\Omega_-(t)\cup\Gamma(t)\cup\Omega_+(t),\\
\Omega_-(t)&:=&\{(x',x_n)\,|\,x'\in\TTM^{n-1},\;0<x_n<h+\eta(x',t)\},\\
\Gamma(t)&:=&\{(x',x_n)\,|\,x'\in\TTM^{n-1},x_n=h+\eta(x',t)\},\\
\Omega_+(t)&:=&\{(x',x_n)\,|\,x'\in\TTM^{n-1},\;h+\eta(x',t)<x_n<L\},
\eea
where $h\in(0,L)$ is a fixed reference level and $\eta\in
C(\RRM^{n-1}\times[0,T])$ is such that $h+\eta(x',t)\in(0,L)$ for all $(x',t)\in
\RRM^{n-1}\times[0,T])$ (cf. Fig. \ref{notation}.)

Following \cite{ilshar, iltscat}, we consider the situation in which the upper
phase $\Omega_-(t)$ is saturated by water under hydrodynamic pressure $P$ while
in the lower phase $\Omega_+(t)$ the pores of the medium are filled by a
vapor-air mixture. This mixture is charaterized by its humidity 
function $\nu$ given by
\[\nu=\frac{\rho_v}{\rho_v+\rho_a}\approx\frac{\rho_v}{\rho_a},\]
where $\rho_v$ and $\rho_a$ are the (variable) density of vapor and the 
(constant) density of air, respectively.  
 The temperature of the mixture and its pressure $P_a$ at the interface are 
assumed constant in time and space. The bulk equations are given just by 
Darcy's law with gravity, incompressibility of water, and constant porosity $m$ 
of the medium, and linear vapor diffusion. Boundary conditions on 
$\Gamma(t)$ express the pressure balance and fixed evaporation/condensation 
humidity $\nu^\ast$. A further condition on $\Gamma(t)$ determines its motion 
from the mass flux balance of water in liquid or 
vapor form across the phase boundary. From this we get the 
complete system \cite{iltscat}
\be\label{mbp0}\left.\begin{array}{rcll}
\Delta P&=&0&\mbox{ in $\Omega_-(t)$,}\\
(\partial_t-D\Delta)\nu&=&0&\mbox{ in $\Omega_+(t)$,}\\
\displaystyle\left(1-\frac{\rho_v}{\rho_w}\right)V_n
&=&\displaystyle-\frac{k}{m\mu_w}
\partial_ { n_\Gamma(t)}
(P-\rho_wgz)+D\frac{\rho_a}{\rho_w}\partial_{n_\Gamma(t)}\nu&
\mbox{ on $\Gamma(t)$,}\\
\nu&=&\nu^\ast&\mbox{ on $\Gamma(t)$,}\\
P&=&P_a+P_c&\mbox{ on $\Gamma(t)$,}\\
P&=&P_0&\mbox{ on $\Sigma_-$,}\\
\nu&=&\nu_a&\mbox{ on $\Sigma_+$,}
\end{array}\right\}
\ee
where  $\Sigma_-:=\{(x',0)\,|\,x'\in\TTM^{n-1}\}$,
$\Sigma_+:=\{(x',L)\,|\,x'\in\TTM^{n-1}\}$, $V_n$ is the normal velocity of 
$\Gamma(t)$, taken positive if $\Omega_-(t)$ is expanding, and $n_{\Gamma(t)}$ 
is the unit normal to $\Gamma(t)$, exterior to $\Omega_-(t)$.

The following
additional constants occur:
\begin{itemize}
\item[$m$:] porosity of the medium ($m\in(0,1)$, fraction of free pore space)
\item[$k$:] its permeability to water,
\item[$\mu_w$:] viscosity of water, 
\item[$D$:] diffusivity of vapor, 
\item[$\rho_w$:] density of water,
\item[$g$:] gravity, 
\item[$P_c$:] capillary pressure, 
\item[$P_0$:] hydrodynamic pressure at upper boundary,
\item[$\nu_a$:] humidity at lower boundary.
\end{itemize}
(Observe that $\rho_v$, which is not constant in the bulk, occurs explicitly
only at  on $\Gamma(t)$ where we have $\rho_v=\rho_a\nu\ast={\rm const}$.)
\subsection{Nondimensionalization}

Substituting $\nu-\nu_a\rightarrow\nu$, $P-P_0\rightarrow P$,
we get 
\be\label{mbp1}\left.\begin{array}{rcll}
\Delta P&=&0&\mbox{ in $\Omega_-(t)$,}\\
(\partial_t-D\Delta)\nu&=&0&\mbox{ in $\Omega_+(t)$,}\\
\left(1-\frac{\rho_v}{\rho_w}\right)V_n&=&-\frac{k}{m\mu_w}\partial_{n_{
\Gamma(t) }}
(P-\rho_wgz)+D\frac{\rho_a}{\rho_w}\partial_{n_{\Gamma(t)}}\nu&
\mbox{ on $\Gamma(t)$,}\\
\nu&=&\nu^\ast-\nu_a&\mbox{ on $\Gamma(t)$,}\\
P&=&P_a+P_c-P_0&\mbox{ on $\Gamma(t)$,}\\
P&=&0&\mbox{ on $\Sigma_-$,}\\
\nu&=&0&\mbox{ on $\Sigma_+$.}
\end{array}\right\}
\ee
We choose $L$ as characteristic length. The characteristic time $T$ and mass $M$ 
are defined in view of (\ref{mbp1})$_3$, (\ref{mbp1})$_5$, from a 
characteristic pressure and velocity
\[\frac{M}{LT^2}=\rho_wgL,\quad \frac{L}{T}=\frac{k\rho_wg}{m\mu_w}.\]
This yields the dimensionless formulation 
\be\label{mbp-nondim}\left.\begin{array}{rcll}
\Delta P&=&0
&\mbox{ in $\Omega_-(t)$,}\\
(\partial_t-\gamma\Delta)\nu&=&0
&\mbox{ in $\Omega_+(t)$,}\\
\left(1-\frac{\rho_v}{\rho_w}\right)V_n&=&-\partial_{n_{\Gamma(t)}}P 
+n_\gamma\cdot e_z
+\frac{\beta}{\nu^\ast-\nu_a}\partial_{n_{\Gamma(t)}}\nu
&\mbox{ on $\Gamma(t)$,}\\
\nu&=&\nu^\ast-\nu_a
&\mbox{ on $\Gamma(t)$,}\\
P&=&\alpha
&\mbox{ on $\Gamma(t)$,}\\
P&=&0
&\mbox{ on $\Sigma_-$,}\\
\nu&=&0
&\mbox{ on $\Sigma_+$}
\end{array}\right\}
\ee
with the dimensionless numbers \cite{iltscat}
\[\alpha=\frac{P_a+P_c-P_0}{\rho_wgL},\quad
\beta=\frac{D\rho_a(\nu^\ast-\nu_a)m\mu_w}{k\rho_w^2gL},\quad
\gamma=\frac{DT}{L^2}=\frac{Dm\mu_w}{k\rho_wgL},\]
see \cite{iltscat} for a physical interpretation of $\alpha$ and $\beta$. 
Denoting the scaled function $\eta$ and the scaled region $\Omega$ by the same symbols, the moving interface is now described by  
\begin{equation}\label{mb}
\Gamma(t)=\{(x',H+\eta(x',t))\} \subset \TTM^{n-1} \times (0,1), \qquad H=h/L\in(0,1),
\end{equation}   
the moving domains are 
\[
\Omega_-(t)=\{(x',x_n)\,|\;0<x_n<H+\eta(x',t)\}, \quad \Omega_+(t)=\{(x',x_n)\,|\;H+\eta(x',t) < x_n < 1\},
\]
enclosed by the fixed hyperplanes 
\[
\Sigma_-=\{(x',0)\}, \qquad \Sigma_+=\{(x',1)\}
\]
and subject to 
\[
\Omega_+(t) \cup \Omega_-(t) \cup \Gamma(t) = \Omega = \TTM^{n-1} \times (0,1).
\]

Normalizing again $\frac{P}{\alpha}\rightarrow P$, $\frac{\nu}{\nu^\ast-\nu_a}\rightarrow \nu$ and scaling once more in the time variable finally leads to the system  
\be\label{mbp4}\left.\begin{array}{rcll}
\Delta P&=&0
&\mbox{ in $\Omega_-(t)$,}\\
(\partial_t-\Delta)\nu&=&0
&\mbox{ in $\Omega_+(t)$,}\\
\mu V_n&=&-\alpha\partial_{n_{\Gamma(t)}}P +
 n_{\Gamma(t)}\cdot e_z
+\beta\partial_{n_{\Gamma(t)}}\nu
&\mbox{ on $\Gamma(t)$,}\\
\nu&=&1
&\mbox{ on $\Gamma(t)$,}\\
P&=&1
&\mbox{ on $\Gamma(t)$,}\\
P&=&0
&\mbox{ on $\Sigma_-$,}\\
\nu&=&0
&\mbox{ on $\Sigma_+$,}
\end{array}\right\}
\ee
which we complement by the initial conditions
\[\eta(\cdot,0)=\eta_0\mbox{ on $\TTM^{n-1}$},\quad\nu(\cdot,0)=\nu_0\mbox{ in $\Omega_+(0)$}.\]
Moreover, we impose the well-posedness condition  
\begin{equation}\label{sign}
\mu \partial_{n_{\Gamma(0)}} \big{[} \beta \nu_0 + \alpha P|_{t=0} \big{]} \leq 
-\omega_0 < 0 \mbox{ on } \Gamma(0).
\end{equation}
Observe that this is in fact a demand on $\eta_0$ and $\nu_0$ only. For later use we collect the facts that 
\[\partial_{n_{\Gamma(t)}}=\frac{1}{\sqrt{1+|\nabla_x \eta|^2}}
(-\nabla_x\eta\cdot\nabla_x+\partial_z),
\quad 
V_n=\frac{\partial_t \eta}{\sqrt{1+|\nabla_x\eta|^2}}\]
and that (\ref{mbp4})$_3$ takes the form 
\[\mu \partial_t\eta=(-\nabla_x\eta\cdot\nabla_x+\partial_z)(-\alpha 
P+\beta\nu)+1
\quad\mbox{ on $\Gamma(t)$.}\]

\section{Transformation} \label{sectrans} 
Following a standard approach we aim to transform system (\ref{mbp4}) to a fixed reference geometry. Oriented at \cite{solfr} we define 
\begin{eqnarray*}
\Omega_-:=\{(x',x_n)\,|\,x_n\in (0,H)\},\quad \Omega_+:=\{(x',x_n)\,|\,x_n\in (H,1)\}, \quad \Gamma:=\{(x',H)\}, \\
%\Omega:=\Omega_-\cup\Omega_+\cup\Gamma, \qquad \qquad \qquad \qquad \qquad \qquad \qquad  
\end{eqnarray*}
and consider continuous functions $\hat\phi:\Omega\times[0,T]\longrightarrow\RRM$ such that
\be\label{phihat}
\hat\phi(\cdot,t)=0\quad\mbox{ on $\Sigma_\pm$},\quad 
\hat\phi((z',H),t)=\eta(z',t),\quad t\in[0,T],
\ee 
 $\hat\phi|_{\Omega_{\pm}}$ and $\hat\phi|_{\Gamma}$ are sufficiently smooth, 
and
\be\label{ddiff}
z\mapsto(z',z_n+\hat\phi(z,t))\in{\rm Diff}(\Omega_{\pm},\Omega_{\pm}(t)),\quad t\in[0,T].
\ee
 (The function $\hat\phi$ we are going to construct in the following will 
satisfy these demands, see Lemma \ref{sigma}, Eqns. (\ref{defphivp}), (\ref{consphi}) and Theorem \ref{mainthrm} below.) Denote the inverse of the mapping (\ref{ddiff}) by $Z_{\hat\phi}(\cdot,t)$ and define 
\begin{eqnarray} \label{Nue} \hat P(z,t)&=&P(Z_{\hat\phi}^{-1}(z,t),t),\quad z\in\Omega_-,\\
\label{P} \hat\nu(z,t)&=&\nu(Z_{\hat\phi}^{-1}(z,t),t),\quad z\in\Omega_+.
\end{eqnarray}
Then system (\ref{mbp4}) is transformed to 
\be\label{mbp5}\left.\begin{array}{rcll}
\displaystyle\cla_\hphi\hat P-\frac{\hat 
P_{z_n}}{1+\hphi_{z_n}}\cla_\hphi\hphi&=&0
&\mbox{ in $\Omega_-\times J$,}\\
\displaystyle L_\hphi\hat \nu-\frac{\hat 
\nu_{z_n}}{1+\hphi_{z_n}}L_\hphi\hphi&=&0
&\mbox{ in $\Omega_+\times J$,}\\
\mu\partial_t\hphi&=&\displaystyle(1+|\nabla'\hphi|^2)\left(\frac{-\alpha\hat 
P_{z_n}}{1+\hphi^-_{z_n}}+\frac{\beta\hat\nu_{z_n}}{1+\hphi^+_{z_n}}\right)+1
&\mbox{ on $\Gamma\times J$,}\\
\hat\nu&=&1
&\mbox{ on $\Gamma\times J$,}\\
\hat P&=&1
&\mbox{ on $\Gamma\times J$,}\\
\hat P&=&0
&\mbox{ on $\Sigma_-\times J$,}\\
\hat\nu&=&0
&\mbox{ on $\Sigma_+\times J$,}\\
\hat\nu(\cdot,0)&=&\hat\nu_0:=\nu_0\circ Z^{-1}_{\hphi}(\cdot,0)
&\mbox{ in $\Omega_+$,}\\
\hphi(\cdot,0)&=&\eta_0&\mbox{ on $\Gamma$,}
\end{array}\right\}
\ee
where 
\bea
L_\hphi u&:=&\partial_t u-\cla_\hphi u,\\
\cla_\hphi u&:=&
\sum_{i=1}^{n-1} u_{z_iz_i}-\vec a(\nabla\hphi)\cdot\nabla u_{z_n},\\
\vec a(\nabla\hphi)&:=&\left(\frac{2\nabla'\hphi}{1+\hphi_{z_n}},
-\frac{1+|\nabla'\hphi|^2}{(1+\hphi_{z_n})^2}\right)
\eea
(Observe that we will assume continuity but not differentiability of $\hphi$ across 
$\Gamma$ and therefore have to distinguish one-sided derivatives $\hphi^\pm_{z_n}$ on $\Gamma$).

\vspace{5mm}

Before we can start to discuss system (\ref{mbp5}) we have to introduce some notation and make some general assumptions which we keep fixed afterwards:
 
\vspace{5mm}

Let $\HHM^n_{\pm}:=\{(x',x_n)\,|\,x' \in \RRM^{n-1},\,x_n\gtrless H\}$, $\RRM^n_\pm:=\{(x',x_n)\,|\,x' \in \RRM^{n-1},\,x_n\gtrless 0\}$.

\vspace{5mm}

For $s \geq 0$ and a Banach space $Y$, $M \in \{ \Omega_\pm, \Gamma, \RRM^n_{\pm}, \HHM^n_{\pm}, [0,T] \}$ ($T > 
0$) we denote by $H_p^s(M,Y)$ the Bessel potential space and by $W_p^s(M,Y)$ the $L^p$-based Sobolev space of order $s$. In 
particular, if $s\notin\mathbb{N}$, this fractional-order Sobolev space coincides with 
the Besov space $B_{pp}^s(M,Y)$, and for $s \in \mathbb{N}$ the we have $W_p^s(M,Y) = H_p^s(M,Y)$ (cf. \cite{trfs1}). For the sake of brevity we write $W_p^s(M):=W_p^s(M,\RRM)$.  

Finally, here and in the following we assume that $p>n+3/2+\sqrt{2n+1/4}$ and that  
\[
\eta_0\in W_p^{2-2/p}(\Gamma), \quad \eta_0(x')\in(\gamma-H,1-\gamma-H) \quad \mbox{for some } \gamma>0.
\]
For technical reasons it is convenient to reduce system (\ref{mbp5}) to the case of homogeneous initial data. For this we need the following two lemmas: 
\begin{lemma}\label{ext0}
There is a linear extension operator $T\in\cll(W_p^{2-1/p}(\Gamma\times\RRM),
H^2_p(\Omega\times\RRM))$ with the properties
\[Tg|_{\Gamma\times\RRM}=g,\quad\partial_nTg|_{\Gamma\times\RRM}=0,\qquad
g\in W_p^{2-1/p}(\Gamma\times\RRM).\]
\end{lemma}
{\bf Proof:} Define first $T^-\in\cll(W_p^{2-1/p}(\Gamma\times\RRM),
H^2_p(\Omega_-\times\RRM))$ by setting (for example) $T^-g:=u$ where $u$ solves 
the elliptic fourth order problem 
\[\Delta^2u=0\quad\mbox{in $\Omega_-\times\RRM$},\quad
u|_{\Gamma\times\RRM}=g,\quad u_{z_n}|_{\Gamma\times\RRM}=0,\quad
u|_{\Sigma_-\times\RRM}=0,\quad u_{z_n}|_{\Sigma_-\times\RRM}=0.\]
Then construct $T$ by extension using \cite{trfs1} Theorem 3.3.4.
\qed
\begin{lemma}\label{sigma}
There is a time interval $J=(0,\tau)$ and $\sigma\in H^2_p(\Omega\times J)$ such that 
%\begin{itemize}
%\item 
\be\label{demsigma}
\sigma|_{\Gamma\times\{0\}}=\eta_0, \quad
\sigma|_{\Sigma_\pm\times J}=0.
\ee
%\item 
Moreover, the mapping $\psi(\cdot,t):=[(z',z_n)\mapsto(z',z_n+\sigma(z,t))]$ is for each $t \in \bar J$ a 
diffeomorphism onto its image satisfying 
\[
\inf_{t \in J} \det D\psi(\cdot,t) \geq \delta>0 
\]
The numbers $\tau$, $\|\sigma\|_{H^2_p(\Omega\times J)}$, and $\delta$ depend 
only on $\|\eta_0\|_{W_p^{2-2/p}(\Gamma)}$ and $\gamma$. 
\end{lemma}
{\bf Proof:}
Define $\eta_1\in W_p^{2-1/p}(\Gamma\times\RRM_+)$ as solution $u$ to the 
standard BVP
\[\Delta u=0\quad\mbox{in $\Gamma\times\RRM_+$,}\quad 
u|_{\Gamma\times\{0\}}=\eta_0.\]
Extend $\eta_1$ to $\eta_2\in W_p^{2-1/p}(\Gamma\times\RRM)$ by \cite{trfs1} 
Theorem 2.9.4 and Proposition 2.9.1.2. Let $h:=\eta_2+H\in 
W_p^{2-1/p}(\Gamma\times\RRM)$. Let $T$ be the operator from Lemma \ref{ext0}.
As $Th$ has H\"older continuous derivatives there are 
$\eta^-\in(H-\gamma/4,H)$, 
$\eta^+\in(H,H+\gamma/4)$, $\tau>0$, such that
\be\label{bounds-}
Th(z,t)>\frac{3}{4}\gamma,\quad 
\partial_nTh(z,t)>-\frac{1}{2}\qquad\mbox{for $z_n\in(\eta^-,H)$, 
$t\in[0,\tau]$},
\ee
\be\label{bounds+}
Th(z,t)<1-\frac{3}{4}\gamma,\quad 
\partial_nTh(z,t)>-\frac{1}{2}\qquad\mbox{for 
$z_n\in(H,\eta^+)$, $t\in[0,\tau]$}.
\ee
Let $\chi\in C^\infty_0(0,1)$ be such that $\supp\chi\in(\eta^-,\eta^+)$,
$\chi(y)\equiv 1$ near $y=H$, $\chi'\geq 0$ on $(0,H)$, and $\chi'\leq 0$ on 
$(H,1)$. 

Define $\psi_n$ by 
\[\psi_n(z,t):=\left\{\begin{array}{rl}
\displaystyle
(1-\chi(z_n))\frac{\gamma}{4}z_n+\chi(z_n)(z_n-H+Th(z,t)),
&\ z\in\Omega_-,\\
\displaystyle
(1-\chi(z_n))(1+\frac{\gamma}{4}(z_n-1))
+\chi(z_n)(z_n-H+Th(z,t)),
&\  z\in\Omega_+\\
\end{array}\right.\]
and let 
\[
\psi(z,t) := (z',\psi_n(z,t)). 
\]
Then 
\[\psi(z',0,t)=0,\quad\psi(z',1,t)=1,\quad\psi(z',H,t)=h,\quad
\partial_n \psi(z',H,t)=1,\qquad t \in J,\]
and to prove the lemma it remains to show that $\partial_n\psi_n\geq\delta>0$ on 
$\Omega$.  
This 
is clear for $z_n\in[0,\eta^-]$ and $z_n\in[\eta^+,1]$. For 
$z_n\in[\eta^-,H]$ we recall $\chi'(z_n)\geq 0$, $z_n-H>-\gamma/4$, and 
(\ref{bounds-}) and conclude
\bea
&&\partial_n\psi_n(z,t)=(1-\chi(z_n))\frac{\gamma}{4}+\chi(z_n)(1+\partial_n
Th(z,t))\\
&&+\chi'(z_n)(z_n-H+Th(z,t)-\frac{\gamma}{4}z_n)\\
&\geq&\min(\gamma/4,1/2).
\eea
Similarly, for $z_n\in[H,\eta^+]$, we have $\chi'(z_n)\leq 0$, $z_n-H<\gamma/4$,
and conclude from this and (\ref{bounds+})
\bea
&&\partial_n\psi_n(z,t)=(1-\chi(z_n))\frac{\gamma}{4}+\chi(z_n)(1+\partial_n
Th(z,t))\\
&&+\chi'(z_n)(z_n-H+Th(z,t)-1-\frac{\gamma}{4}(z_n-1))\\
&\geq&\min(\gamma/4,1/2).
\eea
The proof is completed by setting $\sigma(z,t):=\psi_n(z,t)-z_n$. 
\qed
\begin{remark} \label{regsigma}
{\rm For later use we emphasize the fact that for each $\rho \in (0,1)$ we have 
\[
\sigma \in W_p^\rho((0,\tau),W_p^{2-\rho}(\Omega))
\]
cf. Lemma {\rm4.3} in {\rm \cite{dss}}. }
\end{remark}

\vspace{5mm}

Set 
\begin{equation}\label{defphivp}
\quad\hphi:=\sigma+\phi,\quad \hat P:=Q+p,\quad
\hat\nu=V+v,
\end{equation}
where $Q$ satisfies  the (time-) parameter dependent family of elliptic BVPs
\be\label{defP}\left.\begin{array}{rcll}
\displaystyle
\cla_\sigma Q-\frac{Q_{z_n}}{1+\sigma_{z_n}}\cla_\sigma\sigma&=&0
&\mbox{ in $\Omega_-\times J$,}\\
Q&=&1&\mbox{ on $\Gamma\times J$,}\\
Q&=&0&\mbox{ on $\Sigma_-\times J$,}
\end{array}\right\}
\ee
and $V$ satisfies the parabolic IBVP
\be\label{defV}\left.\begin{array}{rcll}
\displaystyle
L_\sigma V -\frac{V_{z_n}}{1+\sigma_{z_n}}L_\sigma\sigma&=&0
&\mbox{ in $\Omega_+\times J$,}\\
V&=&1&\mbox{ on $\Gamma\times J$,}\\
V&=&0&\mbox{ on $\Sigma_+\times J$,}\\
V(\cdot,0)&=&\hat\nu_0 &\mbox{ in $\Omega_+$.}
\end{array}\right\}
\ee
The functions $Q$ and $V$ are easily seen to be well defined:
\begin{lemma}\label{QV}
Let $\theta \in\left(\frac{1}{p}\,\frac{p-1}{p-n}, \frac{1}{2}\left(1-\frac{n+1}{p}\right)\right)$ and assume $\hat\nu_0 \in 
W_p^{2-2/p}(\Omega_+)$. For some $T \in (0,\tau)$  the problems {\rm 
(\ref{defP})} and {\rm(\ref{defV})} possess unique solutions $Q \in L_p(0,T; 
W_p^2(\Omega_-)) \cap W_p^\theta((0,T),W_p^{2-\theta}(\Omega_-))$ and $V \in 
L_p(0,T;W_p^2(\Omega_+)) \cap H_p^1((0,T);L_p(\Omega_+))$.
\end{lemma}
{\bf Proof:} By construction $1 + \sigma_{z_n}$ is invertible in the Banach 
algebra $H_p^1((0,\tau) \times \Omega)$ and by embedding in the algebras 
$W_p^\theta((0,\tau),W_p^{1-\theta}(\Omega_-))$, 
$L_\infty((0,\tau),W_p^{1-\theta}(\Omega_-))$ and $L_\infty((0,\tau) \times 
\Omega_+)$, too. The assertion is a consequence of the regularity of $\sigma$, 
Lemma \ref{elleq}, and (a periodic version of) known parabolic theory 
\cite{lus}. \qed

\vspace{5mm}

We divide equation $\mbox{(\ref{mbp5})}_{(3)}$ by $\mu \neq 0$ and use the same symbols $\alpha$, $\beta$ again instead of $\alpha / \mu$, $\beta / \mu$. This turns Condition (\ref{sign}) into 
\begin{equation}\label{sign2}
\partial_{n_{\Gamma(0)}} \big{[} \beta \nu_0 + \alpha P|_{t=0} \big{]} \leq 
-\omega_0 < 0 \mbox{ on } \Gamma(0).
\end{equation}
In view of Lemma \ref{QV} the problem (\ref{mbp5}) is reduced to finding $\phi,p,v$ from the 
(formal) system 
\be\label{nonlinF}
F(\sigma+\phi,Q+p,V+v)=0,
\ee
where 
\begin{equation} \label{nonlin0} F(\tau,q,w):=\left(
\begin{array}{c}
\displaystyle\cla_\tau q-\frac{q_{z_n}}{1+\tau_{z_n}}\cla_\tau \tau\\
\displaystyle L_\tau w-\frac{w_{z_n}}{1+\tau_{z_n}} L_\tau\tau\\
\displaystyle\left.\left(\partial_t\tau-(1+|\nabla'\tau|^2)
\left(\frac{-\alpha q_{z_n}}{1+\tau^-_{z_n}}+\frac{\beta 
w_{z_n}}{1+\tau^+_{z_n}}\right)\right)\right|_\Gamma-1/\mu\\
q|_\Gamma-1\\
w|_\Gamma-1\\
q|_{\Sigma_-}\\
w|_{\Sigma_+}
\end{array}\right)\end{equation}
complemented by the given initial conditions. Furthermore, our choices for $\phi$ will ensure that $(\phi,p,v)$ vanishes 
at $t=0$ and is therefore small in suitable norms for short times. Terms that 
are quadratic in this triple will therefore be treated as small perturbations 
to the linearized problem.

Rewrite (\ref{nonlinF}) equivalently as
\be\label{nonlinF2}
F'(\sigma,Q,V)[\phi,p,v]=-F(\sigma,Q,V)-R,
\ee
where 
\bea
R&=&\int_0^1(1-s)\frac{d^2}{ds^2}[F((\sigma,Q,V)+s(\phi,p,v))]\,ds\\
&=&\int_0^1(1-s)F''((\sigma,Q,V)+s(\phi,p,v))[(\phi,p,v),(\phi,p,v)]\,ds\\
&=:&(R_1,R_2,R_3,0,0,0,0).
\eea
Observe, in particular, 
\be\label{aest}
\partial_\tau\cla_\sigma w[\phi]=-A(\nabla\sigma)\nabla\phi\cdot\nabla 
w_{z_n},\quad A(p):=D_p\vec a(p)
\ee

From now on, for the sake of convenience we will write $\phi^\pm$ 
for the restrictions of $\phi$ to 
$\Omega_\pm$ and retain the notation $\phi$ for the trace at $\Gamma$. 
The demands $(\ref{phihat})$ and $(\ref{demsigma})$ together with the 
continuity of $\hat\phi$ across $\Gamma$ imply
\be\label{consphi}
\left.\begin{array}{rcll}
\phi^\pm-\phi&=&0&\mbox{ on $\Gamma\times J$,}\\
\phi^\pm&=&0&\mbox{ on $\Sigma_\pm\times J$,}\\
\phi(\cdot,0)&=&0&\mbox{ on $\Gamma$.}
\end{array}\right\}
\ee
Additionally, it will be convenient to write $Q=:U^-$, $V=:U^+$, $p=u^-$, 
$v=u^+$ $\cla_\sigma=:\Lambda_\sigma^-$, $L_\sigma=:\Lambda_\sigma^+$. Then (\ref{nonlinF2}) takes the form

\be\label{nonlinF3}
\left.\begin{array}{rcll}
\tilde L^\pm u^\pm-\frac{U^\pm_{z_n}}{1+\sigma_{z_n}}L^\pm\phi^\pm
&=&K^\pm\phi^\pm+R^\pm(\phi^\pm,u^\pm)&\mbox{ in $\Omega_\pm\times J$,}\\[1ex]
\partial_t\phi-\alpha^-\phi^-_{z_n}-\alpha^+\phi^+_{z_n}&&&\\
+\zeta \cdot \nabla'\phi
-\tilde\alpha^- u^-_{z_n}-\tilde\alpha^+u^+_{z_n}&=&g_0+
R^B(\phi^-\phi^+,\phi,u^-,u^+)&\mbox{ on $\Gamma\times J$},\\[1ex]
u^\pm&=&0&\mbox{ on $\Sigma_\pm\times J$},\\
u^+(\cdot,0)&=&0&\mbox{ in $\Omega_0$,}
\end{array}\right\}
\ee
where

\bea
L^\pm\phi^\pm&:=&\Lambda_\sigma^\pm\phi^\pm
\pm A(\nabla\sigma)\nabla\phi^\pm\cdot\nabla\sigma_{z_n},\\
\tilde L^\pm u^\pm&:=&\Lambda_\sigma^\pm u^\pm 
-\frac{u^\pm_{z_n}}{1+\sigma_{z_n}}\Lambda_\sigma\sigma,\\
K^\pm\phi^\pm&:=&A(\nabla\sigma)\nabla\phi^\pm\cdot\nabla U^\pm_{z_n}
-\frac{U^\pm_{z_n}}{(1+\sigma_{z_n})^2}\phi^\pm_{z_n}\Lambda_\sigma\sigma,
\eea
\begin{align*}
\alpha^-&:=\alpha \frac{U^-_{z_n}(1+|\nabla'\sigma|^2)}{(1+\sigma_{z_n})^2},
& 
\alpha^+&:=-\beta\frac{ U^+_{z_n}(1+|\nabla'\sigma|^2)}{(1+\sigma_{z_n})^2},\\
\tilde\alpha^-&:= -\alpha\frac{1+|\nabla'\sigma|^2}{1+\sigma_{z_n}},
& \tilde\alpha^+&:= \beta\frac{1+|\nabla'\sigma|^2}{1+\sigma_{z_n}},\\
\zeta&:=\frac{2(\alpha U^-_{z_n}-\beta U^+_{z_n})}{1+\sigma_{z_n}}\nabla'\sigma,&&\\
\end{align*}
\[g_0:=1/\mu - \partial_t\sigma+\frac{1+|\nabla'\sigma|^2}{1+\sigma_{z_n}}
(-\alpha U^-_{z_n}+\beta U^+_{z_n}),\]

\bea
&&R^B(\phi^-,\phi^+,\phi,u^-,u^+)\\
&:=&(1+|\nabla'\sigma|^2)
\left(\alpha\left(
\frac{u_{z_n}^-\phi_{z_n}^-}{(1+\sigma_{z_n})(1+\sigma_{z_n}+\phi^-_{z_n})}
-\frac{U^-_{z_n}{\phi^-_{z_n}}^2}{(1+\sigma_{z_n})^2
(1+\sigma_{z_n}+\phi^-_{z_n})}\right)\right.\\
&&\left.-\beta\left(
\frac{u_{z_n}^+\phi_{z_n}^+}{(1+\sigma_{z_n})(1+\sigma_{z_n}+\phi^+_{z_n})}
-\frac{U^+_{z_n}{\phi^+_{z_n}}^2}{(1+\sigma_{z_n})^2
(1+\sigma_{z_n}+\phi^+_{z_n})}\right)\right)\\
&&+2\nabla'\sigma\cdot\nabla'\phi\left(-\alpha\left(
\frac{u_{z_n}^-}{1+\sigma_{z_n}+\phi^-_{z_n}}
-\frac{U_{z_n}^-\phi_{z_n}^-}{(1+\sigma_{z_n})(1+\sigma_{z_n}+\phi^-_{z_n})}
\right)\right.\\
&&\left.+\beta\left(
\frac{u_{z_n}^+}{1+\sigma_{z_n}+\phi^+_{z_n}}
-\frac{U_{z_n}^+\phi_{z_n}^+}{(1+\sigma_{z_n})(1+\sigma_{z_n}+\phi^+_{z_n})}
\right)\right)\\
&&+|\nabla'\phi|^2\left(-\alpha
\frac{U_{z_n}^-+u_{z_n}^-}{1+\sigma_{z_n}+\phi^-_{z_n}}
+\beta\frac{U_{z_n}^++u_{z_n}^+}{1+\sigma_{z_n}+\phi^+_{z_n}}\right),
\eea
and 
\bea
& & R_\pm(\phi^\pm,u^\pm) \\
&:=& A(\nabla \sigma)\nabla \phi^\pm \nabla u^\pm_{z_n} + \vec{b}(\nabla \sigma, \nabla \phi^\pm) (\nabla U^\pm_{z_n} + \nabla u^\pm_{z_n}) \\
& & + \frac{U^\pm_{z_n} + u^\pm_{z_n}}{1+\sigma_{z_n}+\phi^\pm_{z_n}} \big{(} A(\nabla \sigma)\nabla \phi^\pm \nabla \phi^\pm_{z_n} + \vec{b}(\nabla \sigma, \nabla \phi^\pm) (\nabla \sigma_{z_n} + \nabla \phi^\pm_{z_n}) \big{)} \\
& & + \frac{\phi^\pm_{z_n} L^\pm_\sigma \sigma}{(1+\sigma_{z_n})(1+\sigma_{z_n}+\phi^\pm_{z_n})} (u^\pm_{z_n} - \frac{U^\pm_{z_n} \phi^\pm_{z_n}}{1+\sigma_{z_n}}),  
\eea
where 
\[
\vec{b}(\nabla \sigma, \nabla \phi^\pm) := \sum_{i=0}^n \int_0^1 (1-s) \partial_i A(\nabla \sigma + s \nabla \phi^\pm) \partial_i \phi^\pm \cdot \nabla \phi^\pm \; ds.
\]
To simplify the equations we use the remaining freedom of choice for $\phi$ to 
demand 
\be\label{demphi}
L^\pm\phi^\pm=0\quad\mbox{ in $\Omega_\pm\times J$,}\qquad \phi^+(\cdot,0)=0
\quad\mbox{ in $\Omega_+$.}
\ee
Together with (\ref{consphi}) the demands yield a uniquely solvable time
dependent elliptic boundary value problem for $\phi^-$ and a uniquely solvable
initial boundary value problem for $\phi^+$ that determine $\phi^\pm$
completely in terms of $\phi$. In particular, we have $\phi^-(\cdot,0)=0$ in
$\Omega_-$, and therefore $\sigma=\hphi$ at $t=0$. Comparing now
(\ref{mbp5})$_{1,5,6}$  to (\ref{defP}) at $t=0$ shows $\hat
P(\cdot,0)=Q(\cdot,0)$ in $\Omega_-$, and thus also
$p(\cdot,0)=u^-(\cdot,0)=0$.  

\vspace{5mm}

The following theorem is the main result of this paper:
\begin{theorem} \label{mainthrm}
Let $p$ be as specified above, $\theta\in\left(\frac{1}{p}\,\frac{p-1}{p-n}, \frac{1}{2}\left(1-\frac{n+1}{p}\right)\right)$ and assume {\rm(\ref{sign2})}. There is a positive time $T_0 \in (0,\tau)$ such that the nonlinear problem {\rm(\ref{consphi})}, {\rm(\ref{nonlinF3})}, {\rm(\ref{demphi})} possesses a unique solution 
\begin{eqnarray*}
((u^+,\phi^+),(u^-,\phi^-),\phi) & \in & \quad \big{[} L_p(0,T_0;W_p^2(\Omega_+)) \cap H_p^1((0,T_0);L_p(\Omega_+)) \big{]}^2 \\ 
& & \times \; \big{[} L_p(0,T_0; W_p^2(\Omega_-)) \cap W_p^\theta((0,T_0),W_p^{2-\theta}(\Omega_-)) \big{]}^2 \\ 
& & \times \; L_p(0,T_0;W_p^{2-1/p}(\Gamma))\cap H_p^1((0,T_0),W_p^{1-1/p}(\Gamma)) \\
& & \quad \, \cap \; W^{1+\theta}_p((0,T_0),L_p(\Gamma)).
\end{eqnarray*}
\end{theorem}

\begin{remark}
{\rm We briefly sketch how to obtain a solution of problem {\rm(\ref{mbp4})} from Theorem {\rm \ref{mainthrm}}. Define $\tilde{\phi} \in C(\bar\Omega \times [0,T_0])$ by 
\[ \tilde{\phi} = 
\left\{ \begin{array}{lcr}
\phi^-&\mbox{ in } & \Omega_- \times [0,T_0],\\
\phi^+&\mbox{ in } & \Omega_+ \times [0,T_0],\\
\phi&\mbox{ on } & \Gamma \times [0,T_0],
\end{array}\right.
\]
and let 
\[
\hat\phi := \sigma + \tilde\phi, \qquad \eta := \hat\phi|_{\Gamma \times [0,T_0]}.
\]
Defining $P$ and $\nu$ by means of $\hat{P} := U^- + u^-$, $\hat{\nu} := U^+ + u^+$ and {\rm(\ref{Nue})}, {\rm(\ref{P})}, the triple $(P,\nu,\eta)$ is a solution to {\rm(\ref{mbp4})}.}
\end{remark}

\vspace{5mm}

Our proof of Theorem \ref{mainthrm} will rely on a careful study of the linear part of (\ref{nonlinF3}) which is done in the next section.

\section{The linearized problem} \label{seclin} 
\subsection{Function spaces and optimal regularity for the full linearized problem}
Neglecting lower order terms in (\ref{nonlinF3}) that will be absorbed on the right side we consider the linear elliptic-parabolic problem with dynamic boundary condition
\be\label{lin1}
\left.\begin{array}{rcll}
\displaystyle
\tilde L^\pm u^\pm&=&f^\pm&\mbox{ in $\Omega_\pm$,}\\
L^\pm\phi^\pm&=&0&\mbox{ in $\Omega_\pm$,}\\
\phi^\pm-\phi&=&0&\mbox{ on $\Gamma\times J$,}\\
\phi^\pm&=&0&\mbox{ on $\Sigma_\pm\times J$,}\\
u^\pm&=&0&\mbox{ on $\Gamma\cup\Sigma_\pm$},\\
\partial_t\phi-\alpha^-\phi^-_{z_n}-\alpha^+\phi^+_{z_n}+\zeta \cdot \nabla'\phi-\tilde\alpha^- u^-_{z_n}-\tilde\alpha^+u^+_{z_n}&=&g&\mbox{ on $\Gamma$},\\
u^+(\cdot,0)&=&0&\mbox{ in $\Omega_+$,}\\
\phi^+(\cdot,0)&=&0&\mbox{ in $\bar\Omega_+$.}
\end{array}\right\}
\ee
Here and in the following we assume that   
\[\theta,\theta'\in\left(\frac{1}{p}\,\frac{p-1}{p-n}, 
\frac{1}{2}\left(1-\frac{n+1}{p}\right)\right).\] 
With $J := [0,T]$ ($T \in (0,\tau)$) we further define 
\begin{itemize}
\item for $M\in\{\RRM^{n-1},\TTM^{n-1},\Gamma\}$:
\bea 
X^B(M)&:=&X_\theta^B(M)\\
&:=&L_p(J,W_p^{2-1/p}(M))\cap H_p^1(J,W_p^{1-1/p}(M))
\cap W^{1+\theta}_p(J,L_p(M)),\\
Y^B(M)&:=&Y_\theta^B(M):=L_p(J,W_p^{1-1/p}(M))
\cap W^\theta_p(J,L_p(M)),\\
X_\tr^-(M)&:=&X_{\tr,\theta}^-(M):=L_p(J,W_p^{2-1/p}(M))
 \cap W_p^\theta(J,W^{2-\theta-1/p}_p(M)),
\eea
\item for $M\in\{\HHM^n_+,\Omega_+,\RRM^n_+\}$:
\bea
X^+(M)&:=&L_p(J,H_p^2(M)\cap H_p^1(J,L_p(M)),\\
Y^+(M)&:=&L_p(M\times J),\\
Z^+(M)&:=&L_p(J,H_p^{1}(M))\cap W_p^{1/2}(J,L_p(M)),\\
\eea
\item for $M\in\{\HHM^n_-,\Omega_-,\RRM^n_\pm\}$:
\bea
X^-(M)&:=&X_\theta^-(M):=L_p(J,H_p^2(M)) 
\cap W^\theta_p(J,W^{2-\theta}_p(M)),\\
Y^-(M)&:=&Y_\theta^-(M):=L_p(M\times J)\cap 
W^\theta_p(J,W^{-\theta}_p(M)),\\
Z^-(M)&:=&Z_\theta^-(M):=L_p(J,H_p^{1}(\Omega))\cap 
W_p^\theta(J,W_p^{1-\theta}(\Omega)),
\eea
\item further:
\bea X^+_\tr(\Gamma)&:=& L_p(J,W_p^{2-1/p}(\Gamma))
\cap W_p^{1-1/(2p)}(J,L_p(\Gamma)),\\
Z^0(\Omega)&:=&L_p(J,H_p^{1}(\Omega))\cap H_p^1(J,L_p(\Omega)),\\
Z^+_\tr(\Gamma)&:=&L_p(J,W_p^{1-1/p}(\Gamma))\cap
W_p^{(1-1/p)/2}(J,L_p(\Gamma)),\\
Z^-_\tr(\Gamma)&:=&Z^-_{\tr,\theta}(\Gamma):=L_p(J,W_p^{1-1/p}(\Gamma))\cap 
W_p^\theta(J,W_p^{1-1/p-\theta}(\Gamma)).\\
\eea
\end{itemize}
For any space $U$ in this list, let $\Mathring{U}$ be the (closed) subspace of $U$ 
for which the traces $u|_{t=0}$ and $\partial_t u|_{t=0}$ vanish in case these 
traces exist. The main result of this section reads as follows:  
\begin{theorem} \label{mainthrmlin}
Let Let $p>n+3/2+\sqrt{2n+1/4}$, $\theta,\theta'\in\left(\frac{1}{p}\,\frac{p-1}{p-n}, \frac{1}{2}\left(1-\frac{n+1}{p}\right)\right)$ and assume {\rm(\ref{sign2})}. For some $T_0 \in (0,\tau)$ and 
\[
(f^+,f^-,g) \in Y^+(\Omega_+) \times Y_\theta^-(\Omega_-) \times Y_\theta^B(\Gamma) =: \mathcal{Y}
\]
the problem {\rm(\ref{lin1})} possesses for each $T \in (0,T_0]$ a unique solution 
\[
(u^+,\phi^+,u^-,\phi^-,\phi) \in [ X^+(\Omega_+) ]^2 \times X_\theta^-(\Omega_-) \times X_{\theta'}^-(\Omega_-) \times X^B_\theta(\Gamma) =: \mathcal{X}.
\]
We have $(\phi^-,\phi)(0) = 0$ and the estimate 
\[
\Vert ((u^+,\phi^+),(u^-,\phi^-),\phi) \Vert_{\mathcal{X}} \; \leq \; C \; \Vert (f^+,f^-,g) \Vert_{\mathcal{Y}}
\]
is valid with a constant $C>0$ independent of $T \in (0,T_0]$. If $f^-(0) = 0$ then also $u^-(0) = 0$.
\end{theorem}

%\begin{remark}
%One may be surprised that we obtain solutions to the linear timedependent problem (\ref{lin1}) only for small time intervals. This is due to the localization procedure carried out in the proof of Lemma \ref{varcoeff}. A canonical next step of course would be to perform a continuation procedure producing solutions also on large time intervals. Since we are at the end of day interested in short-time existence of solutions to problem (\ref{nonlinF3}) this remains out of this paper.
%\end{remark}

\subsection{Principal symbol and constant coefficient problems}
As a first step in the proof of Theorem \ref{mainthrmlin} we need to study the principal symbol of problem (\ref{lin1b}). Let $\alpha,\beta\in\RRM$, $c\in\RRM^{n-1}$ such that $\alpha+\beta>0$. 

\vspace{5mm}

For $\kappa\in(0,\pi)$, $\delta\in(0,\pi/2)$ define 
$S_\kappa:=\{re^{i\phi}\,|\,r>0,|\phi|<\kappa\}$,
$\Sigma_\delta:=\{re^{i\phi}\,|\,r\in\RRM\setminus\{0\},|\phi-\pi/2|<\delta\}$.
Define $P\in\Hol(S_\kappa\times\Sigma_\delta^{n-1})$ by
\be\label{defPsym}
P(\lambda,z):=\lambda+\alpha|z|_-+\beta\sqrt{\lambda+|z|_-^2}
-c\cdot z,
\qquad |z|_-:=\sqrt{-\sum_{k=1}^{n-1}z_k^2},
\ee
where $\sqrt{w}:=\sqrt{|w|}e^{i\arg(w)/2}$, $\arg(w)\in(-\pi,\pi]$ for 
$w\in\CCM$. 

$P$ will appear as Fourier symbol of an operator on $\RRM^{n-1}\times\RRM^+$. 
As this symbol is not quasihomogeneous, we determine its $\gamma$ -principal 
part $\pi_\gamma P$ for $\gamma\in (0,\infty]$ in the sense of 
\cite{deka} and find
\[\pi_\gamma P(\lambda,z)=\left\{\begin{array}{cl}
(\alpha+\beta)|z|_--c\cdot z&(\gamma<1),\\
\lambda+(\alpha+\beta)|z|_--c\cdot z&(\gamma=1),\\
\lambda&(\gamma>1).
\end{array}\right.\]
The corresponding Newton polygon is trivial: it is the triangle with vertices 
$(0,0)$, $(1,0)$, $(0,1)$. 
\begin{lemma}\label{lell}
The symbol $P$  is N-parabolic, i.e. there are 
$\delta,\eta\in(0,\pi/2)$ such that 
\be
\label{ell}
\pi_\gamma P(\lambda,z)\neq 0,\quad(\lambda,z)\in S_{\pi/2+\eta}
\times\Sigma^{n-1}_\delta,\quad\gamma\in(0,\infty].
\ee
\end{lemma}
{\bf Proof:} This is trivial for $\gamma>1$. Observe that by homogeneity of 
$\pi_\gamma P$ in $(\lambda,z)$ it is sufficient to show (\ref{ell}) under the 
additional assumption $|z|^2:=\sum|z_k|^2=1$. Furthermore, if 
$z_k\in\Sigma_\delta$ then $|\arg(-z_k^2)|<2\delta$, hence also 
$\left|\arg\left(-\sum_{k=1}^{n-1}z_k^2\right)\right|<2\delta$
and $|\arg|z|_-|<\delta$. On the other hand, $\Re(-z_k^2)\geq 
|z_k|^2\cos(2\delta)$ and therefore 
\[\textstyle||z|_-|^2=||z|_-^2|\geq\Re\left(-\sum_{k=1}^{n-1}
z_k^2\right)\geq\cos(2\delta).\]
Consequently, 
\be\label{estrez-}
\Re(|z|_-)\geq\sqrt{\cos(2\delta)}\cos\delta.
\ee
By the Cauchy-Schwarz inequality and the estimate $|\Re z_k|<|z_k|\sin\delta$
we find
\[\textstyle|\Re(c\cdot z)|=\left|\sum_{k=1}^{n-1} c_k\Re z_k\right|
\leq\|c\|_2\left\|\left(\Re 
z_k\right)_{k=1}^{n-1}\right\|_2\leq\|c\|_2\sin\delta|z|=
\|c\|_2\sin\delta.\]
From this and (\ref{estrez-}),
\be\label{estresym}
\Re((\alpha+\beta)|z|_--c\cdot z)\geq \sqrt{\cos(2\delta)}\cos\delta-
\|c\|_2\sin\delta>\mu>0
\ee 
 for small $\delta>0$. This proves (\ref{ell}) for $\gamma<1$. 

Finally, as $||z|_-|^2=\left|\sum_{k=1}^{n-1}z_k^2\right|\leq|z|^2=1$ it is 
easy to see that 
\[|\Im((\alpha+\beta)|z|_--c\cdot z)|\leq|(\alpha+\beta)|z|_--c\cdot z|
\leq \alpha+\beta+\|c\|_2,\]
and together with (\ref{estresym}) this implies $(\alpha+\beta)|z|_--c\cdot z 
\in S_{\pi/2-\eta}$ for a sufficiently small $\eta>0$. However, if $\lambda\in 
S_{\pi/2+\eta}$ then $-\lambda\notin S_{\pi/2-\eta}$, and this implies 
(\ref{ell}) for $\gamma=1$. \qed

\vspace{5mm}

In the next lemmas we identify $\RRM^{n-1}$ with 
$\partial\RRM^n_+=\RRM^{n-1}\times\{0\}$. 
\begin{lemma}\label{lhspace}(Model problem)

For $g\in 
\Mathring{Y}^B_\theta(\RRM^{n-1})$, there is precisely one 
\[(\phi^-,\phi^+,\phi)\in\Mathring{X}^-_{\theta'}(\RRM^n_+)\times\Mathring{X}
^+(\RRM^n_+)\times\Mathring {X}^B_\theta(\RRM^{n-1})\]
such that 
\be\label{hspace}
\left.\begin{array}{rcll}
\Delta\phi^-&=&0&\mbox{\rm in $\RRM^n_+\times J$,}\\
(\partial_t-\Delta)\phi^+&=&0&\mbox{\rm in $\RRM^n_+\times J$,}\\
\phi^\pm-\phi&=&0&\mbox{\rm on $\RRM^{n-1}\times J$,}\\
\partial_t\phi-\alpha\partial_{x_n}\phi^--\beta\partial_{x_n}\phi^+ 
+c\cdot\nabla'\phi&=&g
&\mbox{\rm on $\RRM^{n-1}\times J$.}
\end{array}\right\}
\ee
There is a constant $C=C(K,\theta,\theta')$ 
such that
\[\|\phi^-\|_{X^-_{\theta'}}+\|\phi^+\|_{X^+}+\|\phi\|_{X^B_\theta}
\leq 
C\|g\|_{Y^B_\theta} ),\]
as long as $\max\{|\alpha|,|\beta|,|c|,(\alpha+\beta)^{-1},T\}\leq K$.
\end{lemma}
{\bf Proof:} Extend $g$ from $J$ to $\RRM_+$ (keeping notation) in such a 
way 
that
\[
\|g\|_{\tilde Y^B_\theta(\RRM^{n-1})}\leq C\|g\|_{Y^B_\theta(\RRM^{n-1})},
\]
where the space $\tilde Y_\theta^B(\RRM^{n-1})$ is obtained from 
$Y_\theta^B(\RRM^{n-1})$
by replacing $J$ by $\RRM_+$. (Observe that the existence of 
the extension and, in particular, the independence of $C$ from $T$ are 
nontrivial and depends on the fact that $g|_{t=0}=0$. We refer to \cite{pss}, 
Proposition 6.1 for the 
details.)

Denote by $(\xi,x_n,\lambda)\mapsto \hat 
\phi^\pm(\xi,x_n,\lambda)$ 
the Fourier transform in the variables $x'$ of the Laplace transform in $t$ of 
(extensions of) $\phi^\pm$ 
and denote by  $\hat\phi$, $\hat g$ the Fourier-Laplace 
transforms 
of  $\phi$ and $g$. Then 
\[
\left.\begin{array}{rcll}
(|\xi|^2-\partial^2_{x_n})\hat \phi^-&=&0&\mbox{ in $\RRM^n_+\times \RRM_+$,}\\
(\lambda+|\xi|^2-\partial^2_{x_n}) \hat \phi^+&=&0&\mbox{ in $\RRM^n_+\times 
\RRM_+$,}\\
\hat \phi^\pm -\hat \phi&=&0&\mbox{ on $\RRM^{n-1}\times 
\RRM_+$,}\\
\lambda\hat\phi_1-\alpha \hat\phi^-_{x_n}-\beta 
\hat\phi^+_{x_n}-ic\cdot\xi \hat \phi
&=&\hat g&\mbox{ on $\RRM^{n-1}\times \RRM_+$.}\\
\end{array}\right\}
\]
As we are seeking regular solutions $\phi^\pm$ this implies
\[\hat \phi^-(\xi,x_n,\lambda)=e^{-|\xi|x_n}\hat \phi^-(\xi,0,\lambda),\quad
\hat \phi^+(\xi,x_n,\lambda)=e^{-\sqrt{\lambda+|\xi|^2}x_n}\hat 
\phi^+(\xi,0,\lambda),\]
and on the boundary
\be\label{bdeq}
P(\lambda,i\xi)\hat \phi(\xi,\lambda)=\hat g(\xi,\lambda)
\ee
with $P$ from (\ref{defPsym}). 
Applying \cite{deka}, Corollary 2.65 and Lemma \ref{lell} we find that there 
exists 
$\omega>0$ and \[\phi \in L_{p,\mbox{\tiny 
loc}}(\RRM_+,W_p^{2-1/p}(\RRM^{n-1}))\cap \Mathring{H}^1_{p,\mbox{\tiny 
loc}}(\RRM_+,W_p^{1-1/p}(\RRM^{n-1}))\cap \Mathring{W}^{1+\theta}_{p,\mbox{\tiny 
loc}}(\RRM_+,W_p^{1-2\theta-1/p}(\RRM^{n-1}))\] such that $\phi$ satisfies (\ref{bdeq}),
 and for $\phi_\omega$ given by 
$\phi_{\omega}(t):=e^{-\omega t}\phi(t)$ 
we have $\phi_{\omega}\in \Mathring{\tilde X}^B_\theta(\RRM^{n-1})$ and
\[\|\phi_{\omega}\|_{\tilde X^B_\theta}\leq C\|g\|_{\tilde Y^B_\theta},\]
where the space  $\tilde X^B_\theta(\RRM^{n-1})$ is
obtained from $X^B_\theta(\RRM^{n-1})$  by replacing $J$ by 
$\RRM_+$. 
Restriction to the interval $J$ yields
\[\|\phi\|_{X^B_\theta}\leq C\|g\|_{Y^B_\theta}.\]

Observe that $X^B_\theta(\RRM^{n-1})\hookrightarrow
X^-_{\tr,\theta'}(\RRM^{n-1})$.
We read (\ref{hspace})$_1$, (\ref{hspace})$_3$ as  Dirichlet problems for 
$\phi^-(\cdot,t)$ and obtain by standard results $\phi^-\in
\Mathring{X}^-_{\theta'}(\RRM^n_+)$ 
and
\[\|\phi^-\|_{X^-_{\theta'}}\leq C\|\phi\|_{X^B_\theta}.\]
Similarly, we read (\ref{hspace})$_1$, (\ref{hspace})$_3$ together with the 
demand $\phi^+|_{t=0}=0$ as an initial-boundary value problem for the heat 
operator solved by $\phi^+$. The compatibility condition occurring in this 
problem is satisfied as $h^+|_{t=0}=0$, and so, by standard results,
$\phi^+\in \Mathring{X}^+$ and
\[\|\phi^+\|_{X_\theta^+}\leq C\|\phi\|_{X_\theta^B}.\]
The statements of the lemma follow now from gathering the given estimates. \qed

\vspace{5mm}

Let $A_0=(a_{ij})\in\RRM^{n\times n}$ be a symmetric positive definite matrix 
with minimal and maximal eigenvalues $\lambda_{\min}$, $\lambda_{\max}$. Let 
$\alpha,\beta\in\RRM$ such that $\beta-\alpha>0$, $c\in\RRM^{n-1}$.
\begin{lemma}\label{fixedcoeff} (Constant coefficients, principal part)

For  $g\in 
\Mathring{Y}^B_\theta(\RRM^{n-1})$, there is precisely one 
\[(\phi^-,\phi^+,\phi)\in\Mathring{X}^-_{\theta'}(\RRM^n_-)\times\Mathring{X}
^+(\RRM^n_+)\times\Mathring {X}^B_\theta(\RRM^{n-1})
\]
such that
\[\left.\begin{array}{rcll}
a_{ij}\partial_{ij}\phi^-&=&0&\mbox{\rm in $\RRM^n_-\times J$,}\\
(\partial_t-a_{ij}\partial_{ij})\phi^+&=&0&\mbox{\rm in $\RRM^n_+\times J$,}\\
\phi^\pm-\phi&=&0&\mbox{\rm on $\RRM^{n-1}\times J$,}\\
\partial_t\phi-\alpha\partial_{x_n}\phi^--\beta\partial_{x_n}\phi^+ 
+c\cdot\nabla'\phi&=&g
&\mbox{\rm on $\RRM^{n-1}\times J$.}
\end{array}\right\}\]
There is a constant $C=C(K,\theta,\theta')$ such that
\[\|\phi^-\|_{X^-_{\theta'}}+\|\phi^+\|_{X^+}+\|\phi\|_{X^B_\theta}
\leq 
C\|g\|_{Y^B_\theta})\]
as long as  $\max\{|\alpha|,|\beta|,|c|,(\beta-\alpha)^{-1}, \lambda_{\max}, 
\lambda_{\min}^{-1},T\}\leq K.$
\end{lemma}
{\bf Proof:} (cf. \cite{lus} \S IV.6) There is an 
$M\in\cll_{is}(\RRM^n)$ 
which leaves $\RRM^n_\pm$ and (hence) $\RRM^{n-1}$ invariant  and satisfies 
$M^\top M=A^{-1}$. 
Consequently, substituting
\[\phi^\pm=\tilde\phi^\pm\circ (M|_{\RRM^n_\pm}\times\id),\quad
\phi=\tilde\phi\circ (M|_{\RRM^{n-1}}\times\id)\]
yields
\[\left.\begin{array}{rcll}
\Delta\tilde\phi^-&=&0&\mbox{\rm in $\RRM^n_-\times J$,}\\
(\partial_t-\Delta)\tilde\phi^+&=&0&\mbox{\rm in 
$\Omega_+\times J$,}\\
\partial_t\tilde\phi-\tilde\alpha\partial_{x_n}\tilde\phi^-
-\tilde\beta\partial_{x_n} \tilde\phi^+ 
+\tilde c\cdot\nabla'\tilde\phi&=&g\circ(M|_{\RRM^{n-1}}\times\id)
&\mbox{\rm on $\Gamma\times J$}
\end{array}\right\}\]
with some $\tilde c\in\RRM^{n-1}$
satisfying $|\tilde c|\leq C(K)$.  Furthermore $\tilde\alpha=M_{nn}\alpha$, 
$\tilde\beta=M_{nn}\beta$, where 
$M_{nn}:=e_n^\top Me_n\in[\lambda_{\max}^{-1},\lambda_{\min}^{-1}]$, so that
\[\tilde\beta-\tilde\alpha\geq\lambda_{\max}^{-1}(\beta-\alpha)>0.\]
To transform the problem to $\Omega_+$ we set 
\[\tilde\phi^-(x',x_n,t)=\bar\phi^-(x',-x_n,t)\]
and obtain a system of the form (\ref{hspace}) with 
$\alpha=-\tilde\alpha$, $\beta=\tilde\beta$, $c=\tilde c$, $\phi^-=\bar\phi^-$,
$\phi^+=\tilde\phi^+$. Now the results follow from Lemma \ref{lhspace} and the
invariance of all occurring function spaces under regular linear 
transformations of the spatial variables. \qed

\subsection{Variable coefficient problems}
We extend the result of Lemma \ref{fixedcoeff} to the case of variable coefficients. Let 
$a_{ij},b_i\in Z^0(\Omega)$, $\alpha\in Z^-_{\tr,\theta}(\Gamma)$, $\beta\in Z^+_{\tr}(\Gamma)$, 
$c\in Y_\theta^B(\Gamma)^{n-1}$ with $a_{ij}=a_{ji}$,
\[a_{ij}(x)\xi^i\xi^j\geq \mu|\xi|^2,\quad x\in\Omega,\,\xi\in\RRM^n,\qquad
\beta(x)-\alpha(x)\geq\mu,\quad x\in\Gamma\]
for some $\mu>0$. 

\begin{lemma}\label{varcoeff} (Variable coefficients, zero initial data)

For any sufficiently small $T>0$ and $g\in 
\Mathring{Y}^B_\theta(\Gamma)$, there is precisely one 
\[(\phi^-,\phi^+,\phi)\in\Mathring{X}^-_{\theta'}(\Omega_-)
\times\Mathring{X}^+(\Omega_+)
\times\Mathring{X}^B_\theta(\Gamma)
\]
such that 
\be\label{varprob}
\left.\begin{array}{rcll}
a_{ij}\partial_{ij}\phi^-+b_i\partial_i\phi^-&=&0&\mbox{\rm in 
$\Omega_-\times J$,}\\
(\partial_t-a_{ij}\partial_{ij})\phi^+-b_i\partial_i\phi^+
&=&0&\mbox{\rm in $\Omega_+\times 
J$,}\\
\phi^\pm-\phi&=&0&\mbox{\rm on $\Gamma\times J$,}\\
\phi^\pm&=&0&\mbox{\rm on $\Sigma_\pm\times J$,}\\
\partial_t\phi-\alpha\partial_{x_n}\phi^--\beta\partial_{x_n}\phi^+ 
+c\cdot\nabla'\phi&=&g
&\mbox{\rm on $\Gamma\times J$.}
\end{array}\right\}
\ee
There is a constant $C=C(K,T,\theta,\theta')$ such that
\[\|\phi^-\|_{X^-_{\theta'}(\Omega_-)}+\|\phi^+\|_{X^+(\Omega_+)}+\|\phi^B\|_{
X^B_\theta(\Gamma) }
\leq 
C\|g\|_{Y^B_\theta(\Gamma)},\]
as long as  
$\max\{\|a_{ij}\|_{Z^0},\|b_i\|_{Z^0},\|\alpha\|_{Z^-_{\tr,\theta}(\Gamma)}, 
\|\beta\|_{Z^+_{\tr}(\Gamma)}, \|c\|_{Y^B_\theta(\Gamma)},\mu^{-1}\}\leq K.$
\end{lemma}
{\bf Proof:} To shorten notation we introduce the operators
\be\label{defLpm}
L^-:=a_{ij}\partial_{ij}+b_i\partial_i,\quad 
L^+:=\partial_t-a_{ij}\partial_{ij}-b_i\partial_i
\ee
It follows from standard results and Lemma \ref{elleq} 
 that we have bounded solution operators
\be\label{defSpm}
\cls^-\in\cll(\Mathring{X}^B_\theta(\Gamma),\Mathring{X}^-_{\theta'}
(\Omega_-)),
\quad
\cls^+\in\cll(\Mathring{X}^B_\theta(\Gamma),\Mathring{X}^+(\Omega_+))
\ee
given for $\psi\in\Mathring{X}^B_\theta(\Gamma)$ by $\cls^\pm\psi:=\psi^\pm$, 
where $\psi^\pm$ solves
\[L^\pm\psi^\pm=0\mbox{ in $\Omega_\pm\times J$,}
\qquad
\psi^\pm=\psi\mbox{ on $\Gamma\times J$,}
\qquad
\psi^\pm=0\mbox{ on $\Sigma_\pm\times J$.}
\]
(Observe that $\psi^-$ has zero time trace automatically while for $\psi^+$ 
this is an additional demand from the choice of the spaces.) 
Further, we define 
$\clt\in\cll(\Mathring{X}^B_\theta(\Gamma),\Mathring{Y}^B_\theta(\Gamma))$
by
\[\clt\psi:=\partial_t\psi-\alpha\partial_{x_n}\cls^-\psi
-\beta\partial_{x_n}\cls^+\psi
+c\cdot\nabla'\psi.\]
To prove the lemma it is sufficient to show that $\clt$ is an isomorphism 
for small $T$ and its inverse has a bound depending only on $K$ and $T$. 
This will be done by the construction of a regularizer, i.e. a map 
$R\in\cll(\Mathring{Y}^B_\theta(\Gamma),\Mathring{X}^B_\theta(\Gamma))$ such 
that
\begin{equation}\label{estreg} \|\clt R-I\|_{\cll(\Mathring{Y}^B_\theta(\Gamma))}\leq 1/2,
\qquad
\|R\clt-I\|_{\cll(\Mathring{X}^B_\theta(\Gamma))}\leq 1/2 \end{equation}
(cf. e.g. \cite {lus} Ch. IV.7/9).

\vspace{5mm}

For small $\lambda>0$, let $\Gamma$ be covered by 
finitely many open balls \newline $\omega^{(k)}=\omega^{(k,\lambda)}:=B(\xi^{(k)},\lambda)$, 
$\Omega^{(k)}=\Omega^{(k,\lambda)}:=B(\xi^{(k)},2\lambda)$ in $\TTM^{n-1} \times (0,1)$, 
$\xi^{(k)}=\xi^{(k,\lambda)}\in\Gamma$,
such that there is an $N_0$ independent of $\lambda$ such that for all $k_0$ 
there are at most $N_0$ balls $\Omega^{(l)}$ with 
$\Omega^{(l)}\cap\Omega^{(k_0)}\neq\varnothing$.  

\vspace{5mm}

Let $\zeta^{(k)}$ be smooth functions with 
$\supp\zeta^{(k)}\subset\Omega^{(k)}$,
$\zeta^{(k)}(x)\in[0,1]$ for all $x\in\Omega$, $\zeta^{(k)}\equiv 1$ on 
$\omega^{(k)}$, $|\partial^\alpha\zeta^{(k)}(x)|\leq 
C_{|\alpha|}\lambda^{-|\alpha|}$. Define additionally $\eta^{(k)}$ by
\[\eta^{(k)}(x):=\frac{\zeta^{(k)}(x)}{\sum_j(\zeta^{(j)}(x))^2},\]
so that
\[\sum_k\eta^{(k)}\zeta^{(k)}\equiv 1.\]
Define $R$ by 
\[Rg:=\sum_k\eta^{(k)}w_k\]
where $g\in\Mathring{Y}^B_\theta(\Gamma)$ and $(w^\pm_k,w_k)$ is the 
solution of the constant-coefficient problem \newline ($\tilde \Gamma := \RRM^{n-1} \times \{ H \}$)
\[\left.\begin{array}{rcll}
L_{0,k}^\pm w_k^\pm&=&0&\mbox{ in $\HHM^n_\pm\times J$,}\\
w^\pm_k-w_k&=&0&\mbox{ on $\tilde \Gamma\times J$,}\\
w^+(\cdot,0)&=&0&\mbox{ on $\tilde \Gamma$,}\\
\clt_{0,k}w_k:=\partial_tw_k-\alpha_{0,k}\partial_{x_n}w_k^-
-\beta_{0,k}\partial_{x_n}w_k^+
+c_{0,k}\nabla'\cdot w_k&=&\zeta^{(k)}g&\mbox{ on $\tilde \Gamma\times J$,}
\end{array}\right\}\]
\[L^-_{0,k}:=a_{ij}(\xi^{(k)})\partial_{ij},\quad
L^+_{0,k}:=\partial_t-a_{ij}(\xi^{(k)})\partial_{ij},\quad
(\alpha_{0,k},\beta_{0,k},c_{0,k}):=(\alpha,\beta,c)(\xi^{(k)}).\]

(Here and in the sequel we identify functions supported in $\Omega^{(k)}$ with 
compactly supported functions on $\tilde\Gamma$.) 

Existence, uniqueness, and estimates for the solution of these problems are 
given in Lemma \ref{fixedcoeff}. Observe, in particular, that $\clt_{0,k}$  is
invertible and 
\[Rg=\sum_k\eta^{(k)}\clt_{0,k}^{-1}(\zeta^{(k)}g).\]
(It is indeed an easy consequence  of Remark 
\ref{suml2}, Lemma \ref{suml4} that the mapping $R\in\cll(\Mathring{Y}^B_\theta(\Gamma),\Mathring{X}^B_\theta(\Gamma))$
is well defined.) 

For later use we note that we have the estimate   
\begin{eqnarray}\label{estwk}
&&\sum_k \|w_k\|_{X^B_\theta(\tilde\Gamma)}^p+\sum_k\|w^+_k\|_{X^+(\HHM^n_+)}^p
+\sum_k\|w^-_k\|_{X^-_{\theta'}(\HHM^n_-)}^p \nonumber\\
&\leq& C \sum_k \|\zeta^{(k)}g\|^p_{Y^B_\theta(\tilde\Gamma)} \leq C \sum_k \|\zeta^{(k)}g\|^p_{Y^B_\theta(\Gamma)} \nonumber\\
&\leq& C(1+\lambda^{-p-n+1} T^{\delta})\|g\|_{Y^B_\theta(\Gamma)}
\end{eqnarray}
for some $\delta > 0$ by Lemmas \ref{fixedcoeff}, \ref{suml4}. 

For an operator $P$ and a function $\phi$,
by $[\phi,P]$ we denote the commutator \newline 
$u\mapsto\phi(Pu)-P(\phi u)$. Choosing smooth cut-off functions $\chi^{(k)} \in \mathcal{D}(\RRM^n)$ such that $\chi^{(k)} \equiv 1$ on $\supp \eta^{(k)}$ and letting $(\tilde{w}_k^\pm,\tilde{w}_k) := \chi^{(k)} (w_k^\pm,w_k)$ we have 
\begin{eqnarray} \label{estreg1}
(\clt R-I)g&=&\sum_k\clt(\eta^{(k)} w_k)-\eta^{(k)}\zeta^{(k)}g \nonumber \\
& = & \sum_k\eta^{(k)}(\clt \tilde w_k-\zeta^{(k)}g)-\sum_k[\eta^{(k)},\clt] \tilde w_k  \\
& = & \sum_k\eta^{(k)}(\clt-\clt_{0,k})\tilde w_k-\sum_k[\eta^{(k)},\clt]\tilde w_k - \sum_k \eta^{(k)} [\chi^{(k)},\clt_{0,k}]w_k. \nonumber
\end{eqnarray}
Thus, in view of (\ref{estreg}), we have to estimate the terms  
\begin{itemize}
\item[1.] $\sum_k\eta^{(k)}(\clt-\clt_{0,k})\tilde w_k$, 
\item[2.] $\sum_k[\eta^{(k)},\clt]\tilde w_k$, and
\item[3.] $\sum_k \eta^{(k)} [\chi^{(k)},\clt_{0,k}]w_k$
\end{itemize}
in $Y^B_\theta(\Gamma)$.

\vspace{5mm} 

1: Let
\[\tilde v_k^\pm:=\cls^\pm \tilde w_k,\quad z^\pm_k:=\cls^\pm(\eta^{(k)}w_k).\]
Using this, we rewrite 
\bea
&&\eta^{(k)}(\clt-\clt_{0,k})\tilde w_k\\
&=&\eta^{(k)}\big(\alpha_{0,k}\partial_{x_n}\tilde w_k^-
+\beta_{0,k}\partial_{x_n}\tilde w_k^+-\alpha\partial_{x_n}\tilde v_k^-
 -\beta\partial_{x_n}\tilde v_k^+
 +(c-c_{0,k})\cdot\nabla' \tilde w_k\big)\\
&=&\eta^{(k)}\big(\alpha\partial_{x_n}(\tilde w_k^--\tilde v_k^-)+\beta\partial_{x_n}
(\tilde w_k^+-\tilde v_k^+)+(\alpha_{0,k}-\alpha)\partial_{x_n}\tilde w_k^-
+(\beta_{0,k}-\beta)\partial_{x_n}\tilde w_k^+\\
&&+(c-c_{0,k})\cdot\nabla' \tilde w_k\big)\\
&=&\alpha\partial_{x_n}(\eta^{(k)}(\tilde w_k^--\tilde v_k^-))
+\beta\partial_{x_n}(\eta^{(k)}(\tilde w_k^+-\tilde v_k^+))
+\eta^{(k)}(\alpha_{0,k}-\alpha)\partial_{x_n}\tilde w_k^-\\
&&+\eta^{(k)}(\beta_{0,k}-\beta)\partial_{x_n}\tilde w_k^+
+\eta^{(k)}(c-c_{0,k})\cdot\nabla' \tilde w_k\\
&&+\alpha[\eta^{(k)},\partial_{x_n}](\tilde w_k^--\tilde v_k^-)
+\beta[\eta^{(k)},\partial_{x_n}](\tilde w_k^+-\tilde v_k^+) \\
&=&\alpha\partial_{x_n}(\eta^{(k)}(\tilde w_k^--\tilde v_k^-))
+\beta\partial_{x_n}(\eta^{(k)}(\tilde w_k^+-\tilde v_k^+))
+\eta^{(k)}(\alpha_{0,k}-\alpha)\partial_{x_n}\tilde w_k^-\\
&&+\eta^{(k)}(\beta_{0,k}-\beta)\partial_{x_n}\tilde w_k^+
+\eta^{(k)}(c-c_{0,k})\cdot\nabla' \tilde w_k\\
%&&+\alpha[\eta^{(k)},\partial_{x_n}](\tilde w_k^--\tilde v_k^-)
%+\beta[\eta^{(k)},\partial_{x_n}](\tilde w_k^+-\tilde v_k^+) \\
\eea
since 
\[
[\eta^{(k)},\partial_{x_n}](\tilde w_k^\pm-\tilde v_k^\pm) = - (\partial_{x_n}\eta^{(k)}) (\tilde w_k^\pm-\tilde v_k^\pm) = 0.
\]
As a consequence of Lemma \ref{suml1} and Remark \ref{suml2} it suffices to estimate 
\begin{itemize}
\item[1.1.] $\eta^{(k)}(\tilde w_k^\pm-\tilde v_k^\pm)$  in $X^\pm(\Omega_\pm)$, 
\item[1.2.] $\eta^{(k)}(\alpha-\alpha_{0,k})\partial_{x_n}\tilde w^-_k$,
$\eta^{(k)}(\beta-\beta_{0,k})\partial_{x_n}\tilde w^+_k$ in 
$Y^B_\theta(\Gamma)$, 
\item[1.3.] $\eta^{(k)}(c-c_{0,k})\nabla' \tilde w_k$ in 
$Y^B_\theta(\Gamma)$. 
\end{itemize}

1.1: Observe that the differences $\eta^{(k)}(\tilde w_k^\pm-\tilde v_k^\pm) = \eta^{(k)}(w_k^\pm-\tilde v_k^\pm)$ solve the boundary value problems
\[\left.\begin{array}{rcll}
L^\pm\big(\eta^{(k)}(\tilde w_k^\pm-\tilde v_k^\pm)\big)
&=&(L^\pm-L_{0,k}^\pm)(\eta^{(k)}\tilde w_k^\pm)\\[1mm]
&& -[\eta^{(k)},L_{0,k}^\pm]\tilde w_k^\pm
+[\eta^{(k)},L^\pm]\tilde v_k^\pm &\mbox{ in $\Omega_\pm\times J,$}\\
\eta^{(k)}(\tilde w_k^\pm-\tilde v_k^\pm)&=&0&\mbox{ on $(\Gamma\cup\Sigma_\pm)\times J$.}
\end{array}\right\}\]

Hence we have to consider 

\begin{itemize}
\item[1.1.1.] $\eta^{(k)}(L^\pm-L^\pm_{0,k})\tilde w_k^\pm$, 
\item[1.1.2.] $[\eta^{(k)},L^\pm]\tilde w_k^\pm$, $[\eta^{(k)},L_{0,k}^\pm]\tilde w_k^\pm$, and 
\item[1.1.3.] $[\eta^{(k)},L^\pm]\tilde v_k^\pm$
\end{itemize}
in $Y^\pm(\Omega_\pm)$.

1.1.1: As the $a_{ij}$ are H\"older continuous with an exponent 
$\kappa>0$ we have
\bea
\|\eta^{(k)}(a_{ij}-a_{ij}(\xi^{k}))\partial_{ij}\tilde w_k^{\pm}\|_{
L^p(\Omega\times J)}
&\leq& 
C\|\eta^{(k)}(a_{ij}-a_{ij}(\xi^{k}))\|_\infty\|\partial_{ij}\tilde w_k^{\pm}\|_{
L^p(\Omega\times 
J)}\\
&\leq& C\lambda^\kappa\|\tilde w_k^\pm\|_{X^\pm(\Omega_\pm)} \leq C\lambda^\kappa\| w_k^\pm\|_{X^\pm(\HHM^{n-1}_\pm)}.
\eea
Fix $\theta_{1,2,3}$ such that 
$\theta<\theta_1<\theta_2<\theta_3<\frac{1}{2}\left(1-\frac{n+1}{p}\right)$.
Using the multiplication property (cf. Lemmas \ref{mult2}, \ref{mult3})
\[BUC^{\theta_1}(J,BUC^{\theta_2}(\Omega_-))\cdot
W_p^\theta(J,W_p^{-\theta}(\Omega_-))\hookrightarrow
W_p^\theta(J,W_p^{-\theta}(\Omega_-))\]
and Lemma \ref{remain} we get 
\bea
&&\|\eta^{(k)}(a_{ij}-a_{ij}(\xi^{k}))\partial_{ij}\tilde w_k^-\|_{
W_p^\theta(J,W_p^{-\theta}(\Omega_-))}\\
&\leq& C
\|\eta^{(k)}\|_{BUC^{\theta_2}(\Omega_-)}\|a_{ij}-a_{ij}(\xi^{k})\|_{
BUC^{\theta_1}(J,BUC^{\theta_2}(\Omega_-\cap\Omega^{(k)}))}\|\partial_{ij}
\tilde w_k^-\|_ {
W_p^\theta(J,W_p^{-\theta}(\Omega_-))}\\
&\leq &C\lambda^{-\theta_2}(T^{\theta_2-\theta_1}+\lambda^{2\theta_3-\theta_2})
 \|a_{ij}\|_{Z^0(\Omega)}\|\tilde w_k^-\|_{X^-(\Omega_-)} \\
&\leq &C\lambda^{-\theta_2}(T^{\theta_2-\theta_1}+\lambda^{2\theta_3-\theta_2})
 \|a_{ij}\|_{Z^0(\Omega)}\| w_k^-\|_{X^-(\HHM_-)}
\eea

1.1.2: Note that formally 
\bea
&&-[\etk,L^-]u=[\etk,L^+]u\\
&=&a_{ij}(\partial_{ij}\etk u+\partial_i\etk \partial_j u+
\partial_j\etk \partial_i u)+b_i\partial_i\etk u,
\eea 
so
\bea
\|[\etk,L^\pm]\tilde w_k^\pm\|_{L_p(\Omega_\pm\times J)}
&\leq& 
C(\|a\|_\infty+\|b\|_\infty)\|\etk\|_{W^2_\infty(\Omega_\pm\times J)}
\|\tilde w_k^\pm\|_{L_p(J,H^1_p(\Omega_\pm))}\\
&\leq& C\lambda^{-2}T^\delta\|\tilde w_k^\pm\|_{X^\pm(\Omega_\pm)}
\leq C\lambda^{-2}T^\delta\|\tilde w_k\|_{X^B(\Gamma)} \\
&\leq& C\lambda^{-2}T^\delta\| w_k\|_{X^B(\tilde \Gamma)}
\eea
for some $\delta > 0$. In the $+$\ -case, this is what we need to show. In the $-$\ -case, we 
additionally use product estimates parallel to those derived in Lemma \ref{mult2} 
%the product estimate {\bf (Check!)}
%\[\|uv\|_{W_p^\theta(J,L_p(\Omega_-))}
%\leq C(\|u\|_\infty\|v\|_{W_p^\theta(J,L_p(\Omega_-))}
%+\|u\|_{W_p^\theta(J,L_p(\Omega_-))}\|v\|_\infty)\]
to get
\bea
&&\|[\etk,L^-]\tilde w_k^-\|_{W_p^\theta(J,W_p^{-\theta}(\Omega_-))}\\
&\leq &C\|[\etk,L^-]\tilde w_k^-\|_{W_p^\theta(J,L_p(\Omega_-))} \\
&\leq& C\lambda^{-2}\big(\|\tilde w_k^-\|_{W_p^\theta(J,H^1_p(\Omega_-))}
+\|\tilde w_k^-\|_{L_\infty(J,W_\infty^1(\Omega_-))}\big)\\
&\leq&C\lambda^{-2}T^\delta\|\tilde w_k^-\|_{X^-_{\theta'}(\Omega_-)} \leq C\lambda^{-2}T^\delta\| w_k^-\|_{X^-_{\theta'}(\HHM^{n-1}_-)} \\
%& \leq & C\lambda^{-2}T^\delta\|\tilde w_k^-\|_{X^-_{\theta'}(\HHM^{n-1}_-)} \\
%&\leq& C \lambda^{-2}T^\delta\|g\|_{Y^B_\theta(\Gamma)}
\eea
for some $\delta > 0$. The terms $[\etk,L_{0,k}^\pm]w_k^\pm$ are treated in the same way.

1.1.3: We use arguments parallel to 1.1.2, using additionally
\[\|\tilde v_k^-\|_{X_{\theta'}^-(\Omega_-)}+\|\tilde v_k^+\|_{X^+(\Omega_+)}\leq C
\|\tilde w_k\|_{X^B_\theta(\Gamma)} \leq C
\|w_k\|_{X^B_\theta(\tilde \Gamma)} \]
by standard parabolic theory and Lemma \ref{elleq}. \\

1.2: Using 
\[\|\partial_{x_n} \tilde w_k^\pm\|_{Y^B(\Gamma)}\leq 
C\|\tilde w_k^\pm\|_{X^\pm(\Omega_\pm)},\]
these terms can be estimated in the same way as the following 

1.3: In general we have by Lemmas \ref{mult2}, \ref{mult3} and \ref{alg} 
\be\label{prod-inf}
\|uv\|_{Y^B(\Gamma)}\leq
C(\|u\|_{Y^B(\Gamma)}\|v\|_\infty+\|u\|_\infty\|v\|_{Y^B(\Gamma)}) \qquad u, v \in\Mathring{Y}^B(\Gamma),
\ee 
and 
\be\label{inf-small}
\|u\|_\infty\leq C T^\delta\|u\|_{Y^B(\Gamma)},\qquad u\in
\Mathring{Y}^B(\Gamma)
\ee
for some $\delta>0$ (cf. also Remark \ref{rem2}). As $c$ is H\"older continuous (in space and time) with exponent 
$\kappa>0$, 
\bea
&&\|\eta^{(k)}(c-c_{0,k})\nabla' \tilde w_k\|_{Y^B(\Gamma)}\\
&\leq& C\big(
\|\eta^{(k)}(c-c_{0,k})\|_\infty\|\tilde w_k\|_{X^B(\Gamma)} + 
\|\eta^{(k)}(c-c_{0,k})\|_{Y^B(\Gamma)}
\|\nabla' \tilde w_k\|_\infty\big)\\
&\leq& C\big(\lambda^\kappa\|\tilde w_k\|_{X^B(\Gamma)}+T^\delta
\|\nabla' \tilde w_k\|_{Y^B(\Gamma)}\big)\leq C(\lambda^\kappa+T^\delta)
\|\tilde w_k\|_{X^B(\Gamma)} \\
&\leq&  C(\lambda^\kappa+T^\delta)
\|w_k\|_{X^B(\tilde \Gamma)}
\eea
for some $\delta > 0$. 

\vspace{5mm}

2: We have 
\[
\Vert \sum_k [\eta^{(k)},\clt]\tilde w_k \Vert_{Y^B_\theta(\Gamma)}^p \leq C(\lambda) \max_{k} \Vert [\eta^{(k)},\clt]\tilde w_k \Vert_{Y^B_\theta(\Gamma)}^p,
\]
\bea
[\eta^{(k)},\clt]\tilde w_k
&=&-\alpha[\eta^{(k)},\partial_{x_n}]\tilde v_k^--\beta[\eta^{(k)},\partial_{x_n}]\tilde v_k^+
\\
&&-\alpha\partial_{x_n}(z_k^--\eta^{(k)}\tilde v_k^-)
-\beta\partial_{x_n}(z_k^+-\eta^{(k)}\tilde v_k^+)
+[\eta^{(k)},c\cdot\nabla']\tilde w_k.
\eea
More explicitly, we get by calculating the commutators and using the definition
of $\tilde v_k$ 
\bea
%[\eta^{(k)},\partial_{x_n}](\tilde w_k^\pm-\tilde v_k^\pm)&=&0,\\
-[\eta^{(k)},\partial_{x_n}]\tilde v_k^\pm&=&\partial_{x_n}\eta^{(k)}\tilde w_k,\\
{}[\eta^{(k)},c\cdot\nabla']\tilde w_k&=&-c\cdot(\nabla'\eta^{(k)})\tilde w_k,
\eea
so we have to consider 
\begin{itemize}
\item[2.1.] $z_k^\pm-\eta^{(k)}\tilde v_k^\pm$ in $X^\pm(\Omega_\pm)$,
\item[2.2.] $c\cdot(\nabla'\eta^{(k)})\tilde w_k$ in $Y^B_\theta(\Gamma)$,
\item[2.3.] $\partial_{x_n}\eta^{(k)}\tilde w_k$ in  $Y^B_\theta(\Gamma)$.
\end{itemize}

2.1: The differences $z_k^\pm-\eta^{(k)}\tilde v_k^\pm$ are solutions to 
\[\left.\begin{array}{rcll}
L^\pm\big(z_k^\pm-\eta^{(k)}\tilde v_k^\pm\big)
&=&[\eta^{(k)},L^\pm]\tilde v_k^\pm &\mbox{ in $\Omega_\pm\times J,$}\\
z_k^\pm-\eta^{(k)}\tilde v_k^\pm&=&0&\mbox{ on $(\Gamma\cup\Sigma_\pm)\times J.$}
\end{array}\right\}\]
Therefore, 2.1 can be estimated parallel to 1.1.2, 1.1.3.

2.2: We have 
\bea
&&\|c\cdot(\nabla'\eta^{(k)})\tilde w_k\|_{Y^B(\Gamma)}\\
&\leq& \|c\cdot(\nabla'\eta^{(k)})\|_{Y^B(\Gamma)}\|\tilde w_k\|_{Y^B(\Gamma)}
\leq C\lambda^{1/p-2}\|\tilde w_k\|_{Y^B(\Gamma)}\leq C 
\lambda^{1/p-2}T^\delta
\|\tilde w_k\|_{X^B(\Gamma)} \\
& \leq & C 
\lambda^{1/p-2}T^\delta
\|w_k\|_{X^B(\tilde \Gamma)}
\eea
for some $\delta > 0$.  

2.3: This is handled in the same fashion as 2.2. 

\vspace{5mm} 

3: Observe that 
\[
[\chi^{(k)},\mathcal{T}_{0,k}]w_k = w_k (c_{0,k} \cdot \nabla' - (\alpha_{0,k} + \beta_{0,k}) \partial_{x_n} ) \chi^{(k)}. 
\]
Using 
\begin{eqnarray*}
& & \Vert \sum_k \eta^{(k)} (w_k (c_{0,k} \cdot \nabla' - (\alpha_{0,k} + \beta_{0,k}) \partial_{x_n} ) \chi^{(k)}) \Vert_{Y^B_\theta( \Gamma)}^p \\  
& \leq & C(\lambda) \max_{k} \Vert w_k (c_{0,k} \cdot \nabla' - (\alpha_{0,k} + \beta_{0,k}) \partial_{x_n} ) \chi^{(k)} \Vert_{Y^B_\theta(\tilde \Gamma)}^p
\end{eqnarray*}
as well as Lemmas \ref{mult2}, \ref{mult3}, \ref{suml1}, \ref{suml4} and Remark \ref{rem1} it is easily verified that 
\begin{eqnarray*}
& & \Vert \sum_k w_k (c_{0,k} \cdot \nabla' - (\alpha_{0,k} + \beta_{0,k}) \partial_{x_n} ) \chi^{(k)} \Vert_{Y^B_\theta(\tilde \Gamma)}^p \\  
& \leq & C(\lambda) T^\delta  \max_{k} \|w_k\|_{X^B_\theta(\tilde \Gamma)}^p  \leq C(\lambda) T^\delta  \max_{k} \|\zeta^{(k)}g\|^p_{Y^B_\theta(\tilde \Gamma)} \\
& \leq & C(\lambda) T^\delta  \max_{k} \|g\|^p_{Y^B_\theta(\Gamma)},
\end{eqnarray*}
where the constant $C(\lambda)$ may differ from term to term.

\vspace{5mm}

%Using that $k \sim \lambda^{-n+1}$ and
%\begin{eqnarray*}
%& & \int_{\Omega_k \cap \Gamma} \int_{\Omega_k \cap \Gamma} \frac{|\zeta^{(k)}(x)-\zeta^{(k)}(y)|^p}{|x-y|^{n-1+rp}} \; dx\;dy \\
%& \leq & C \, \lambda^{-p} \, \int_{\mathbb{T}^{n-1}} \int_{\mathbb{T}^{n-1}} |x-y|^{-n+1+p(1-r)} \; dx\;dy \\
%& \leq & C \, \lambda^{-p}
%\end{eqnarray*}
%($r:=1-1/p$) direct calculations show that
%\begin{eqnarray*}
%\sum_k \int_0^T \Vert \zeta^{(k)} g(t) \Vert^p_{L_p(\Gamma)} \, dt & \leq & C(N_0)  \, \Vert g \Vert^p_{L_p(0,T; L_p(\Gamma))}, \\
%\sum_k \int_0^T [ \zeta^{(k)} g(t) ]^p_{p,r,\Omega_k \cap \Gamma} \, dt & \leq & C(N_0)  \, \Vert g \Vert^p_{L_p(0,T;W^r_p(\Gamma)} \\ 
%& & + \, C \, \lambda^{-p-n+1} \, T \, \Vert g \Vert^p_{L_\infty((0,T) \times \Gamma)}, \\
%\sum_k [ \zeta^{(k)} g ]^p_{p,\theta,0,T;L_p(\Gamma)} & \leq & C(N_0) \, [ g ]^p_{p,\theta,0,T;L_p(\Gamma)}.
%\end{eqnarray*}
%
%
%Lemma \ref{fixedcoeff}, (\ref{inf-small}) and these estimates yield 
%\begin{eqnarray}\label{estwk}
%&&\sum_k \|w_k\|_{X^B_\theta(\Gamma)}^p+\sum_k\|w^+_k\|_{X^+(\Omega_+)}^p
%+\sum_k\|w^-_k\|_{X^-_{\theta'}(\Omega_-)}^p \nonumber\\
%&\leq& C \sum_k \|\zeta^{(k)}g\|^p_{Y^B_\theta(\Gamma)}\nonumber\\
%&\leq& C(1+\lambda^{-p-n+1} T^{1+\delta})\|g\|_{Y^B_\theta(\Gamma)}.
%\end{eqnarray}
%
%\vspace{5mm}
%
The first estimate in (\ref{estreg}) now follows from (\ref{estwk}) by choosing first $\lambda > 0$ and then $T > 0$ small enough.

\vspace{5mm}

Reversely, we have for $u\in\Mathring{X}^B_\theta(\Gamma)$ that 
\[(R\clt-I)u=\sum_k\big(\eta^{(k)}T_{0,k}^{-1}(\zeta^{(k)}\clt
u)-\eta^{(k)}\zeta^{(k)}u\big)
=\sum_k\eta^{(k)}T_{0,k}^{-1}\big(\zeta^{(k)}\clt
u-\clt_{0,k}(\zeta^{(k)}u)\big)\]
and 
\[\zeta^{(k)}\clt u-\clt_{0,k}(\zeta^{(k)}u)=
\zeta^{(k)}(\clt-\clt_{0,k})(\chi^{(k)}u)+[\zeta^{(k)},\clt_{0,k}](\chi^{(k)}u)+\zeta^{(k)}[\chi^{(k)},\clt]u.\]
The second estimate in (\ref{estreg}) can be obtained using (\ref{sum12}) in Remark \ref{suml5}, Lemma \ref{fixedcoeff} and arguments parallel to those used to treat the terms $\eta^{(k)}(\clt-\clt_{0,k})\tilde w_k$, $[\chi^{(k)},\clt_{0,k}]w_k$ and $[\eta^{(k)},\clt]\tilde w_k$ above. This proves the lemma. \qed

\begin{lemma} \label{inhomindat} (Inhomogeneous initial data)

For any sufficiently small $T>0$ and $g\in 
Y^B_\theta(\Gamma)$, there is precisely one 
\[(\phi^-,\phi^+,\phi)\in\Mathring{X}^-_{\theta'}(\Omega_-)
\times\Mathring{X}^+(\Omega_+)
\times X^B_\theta(\Gamma)\]
satisfying {\rm(\ref{varprob})}.
There is a constant $C=C(K,T,\theta,\theta')$ such that
\[\|\phi^-\|_{X^-_{\theta'}(\Omega_-)}+\|\phi^+\|_{X^+(\Omega_+)}+\|\phi^B\|_{
X^B_\theta(\Gamma)}
\leq 
C\|g\|_{Y^B_\theta(\Gamma)},\]
as long as  
$\max\{\|a_{ij}\|_{Z^0},\|b_i\|_{Z^0},\|\alpha\|_{Z^-_{\tr,\theta}(\Gamma)}, 
\|\beta\|_{Z^+_{\tr}(\Gamma)}, \|c\|_{Y^B_\theta(\Gamma)},\mu^{-1}\}\leq K.$
\end{lemma}
{\bf Proof:} From \cite{dss}, Theorem 4.5., it follows that there is a
$\tilde\phi=\tilde\phi(g)\in X_\theta^B$ such that 
\[\tilde\phi|_{t=0}=0,\quad \partial_t\tilde\phi|_{t=0}=g,\quad
\|\tilde\phi\|_{X_\theta^B}\leq C\|g\|_{Y_\theta^B}\]
with $C$ independent of $g$. Using  the solution
operators defined in (\ref{defSpm}) we split
\[(\phi^-,\phi^+,\phi)=(\phi_1^-,\phi_1^+,\phi_1)
+(\cls^-\tilde\phi,\cls^-\tilde\phi,\tilde\phi),\]
where $(\phi_1^-,\phi_1^+,\phi_1)\in\Mathring{X}^-_{\theta'}\times\Mathring{X}^+
\times\Mathring{X}^B_\theta(\Gamma)$ satisfies (cf. (\ref{defLpm}))
\be\label{phi1}
\left.\begin{array}{rcll}
L^\pm\phi_1^\pm&=&0&\mbox{ in $\Omega_\pm\times J$,}\\
\phi_1^\pm-\phi_1&=&0&\mbox{ on $\Gamma\times J$,}\\
\phi_1^\pm&=&0&\mbox{ on $\Sigma_\pm\times J$,}\\
\partial_t\phi_1-\alpha\partial_{x_n}\phi_1^--\beta\partial_{x_n}\phi_1^+
+c\cdot\nabla'\phi_1&=&g_1&\mbox{ on $\Gamma\times J$,}
\end{array}\right\}
\ee
with
\[g_1:=g-\partial_t\tilde\phi+\alpha\partial_{x_n}\cls^-\tilde\phi
+\beta\partial_{x_n}\cls^+\tilde\phi-c\cdot\nabla'\tilde\phi\in\Mathring{Y}
_\theta^B(\Gamma).\]
From (\ref{defSpm}) we have
\[\|\cls^-\tilde\phi\|_{X^-_{\theta'}(\Omega_-)}\leq
C\|g\|_{Y^B_\theta(\Gamma)},\quad
\|\cls^+\tilde\phi\|_{X^+}(\Omega_+)\leq C\|g\|_{Y^B_\theta(\Gamma)},\]
and consequently
\[\|g_1\|_{Y^B_\theta(\Gamma)}\leq C \|g\|_{Y^B_\theta(\Gamma)}.\]
The lemma follows from this and the application of Lemma \ref{varcoeff} to 
the system (\ref{phi1}).\qed \vspace{5mm}
{\bf Proof of Theorem \ref{mainthrmlin}:} System (\ref{lin1}) splits into the problems 
\be\label{lin1a}
\left.\begin{array}{rcll}
\displaystyle
\tilde L^\pm u^\pm&=&f^\pm&\mbox{ in $\Omega_\pm$,}\\
%L^\pm\phi^\pm&=&0&\mbox{ in $\Omega_\pm$,}\\
%\phi^\pm-\phi&=&0&\mbox{ on $\Gamma\times J$,}\\
%\phi^\pm&=&0&\mbox{ on $\Sigma_\pm\times J$,}\\
u^\pm&=&0&\mbox{ on $\Gamma\cup\Sigma_\pm$},\\
%\partial_t\phi-\alpha^-\phi^-_{z_n}-\alpha^+\phi^+_{z_n}+c\cdot\nabla'\phi-\tilde\alpha^- u^-_{z_n}-\tilde\alpha^+u^+_{z_n}&=&g&\mbox{ on $\Gamma$},\\
u^+(\cdot,0)&=&0&\mbox{ in $\Omega_+$,}\\
%\phi(\cdot,0)&=&0&\mbox{ on $\Gamma$,}
\end{array}\right\}
\ee
and
\be\label{lin1b}
\left.\begin{array}{rcll}
\displaystyle
%\tilde L^\pm u^\pm&=&f^\pm&\mbox{ in $\Omega_\pm$,}\\
L^\pm\phi^\pm&=&0&\mbox{ in $\Omega_\pm$,}\\
\phi^\pm-\phi&=&0&\mbox{ on $\Gamma\times J$,}\\
\phi^\pm&=&0&\mbox{ on $\Sigma_\pm\times J$,}\\
%u^\pm&=&0&\mbox{ on $\Gamma\cup\Sigma_\pm$},\\
\partial_t\phi-\alpha^-\phi^-_{z_n}-\alpha^+\phi^+_{z_n}+ \zeta \cdot\nabla'\phi&=&g+\tilde\alpha^- u^-_{z_n}+\tilde\alpha^+u^+_{z_n}&\mbox{ on $\Gamma$},\\
%u^+(\cdot,0)&=&0&\mbox{ in $\Omega_+$,}\\
\phi^+(\cdot,0)&=&0&\mbox{ in $\bar\Omega_+$.}
\end{array}\right\}
\ee
Moreover, Condition (\ref{sign2}) implies that 
\[
(\alpha_+ -\alpha_-)(t) = \big{[} (-\beta U^+_{z_n} - \alpha U^-_{z_n}) \frac{ 
1+|\nabla'\sigma|^2}{(1+\sigma_{z_n})^2} \big{]}(t) \geq \omega_1 > 0 \quad 
\mbox{on } \Gamma 
\]
for small $t \geq 0$ by continuity. The assertion follows from first applying Lemma \ref{elleq} and standard parabolic theory  to (\ref{lin1a}) and then applying Lemma \ref{inhomindat} to (\ref{lin1b}). \qed

\section{The nonlinear problem} \label{secnonlin} 
Let $p$ be as specified in Section $3$ and  
$\theta,\theta'\in\left(\frac{1}{p}\,\frac{p-1}{p-n}, 
\frac{1}{2}\left(1-\frac{n+1}{p}\right)\right)$ be such that $\theta' > \theta$. 
Observe that  
\[
\theta + \theta' < 1 - n/p, \quad  \mbox{ hence } \quad 
W_p^{1-\theta'}(\Omega_-) \cdot W_p^{-\theta}(\Omega_-) \hookrightarrow 
W_p^{-\theta}(\Omega_-). 
\]
Recall the definitions of the spaces $\mathcal{X}$ and $\mathcal{Y}$ from 
Theorem \ref{mainthrmlin} , and let $L^\pm$, $\tilde{L}^\pm$, $\alpha^\pm$, 
$\tilde{\alpha}^\pm$, $\zeta$ be the 
operators and coefficient functions introduced in (\ref{nonlinF3}).
Let $\tilde{\mathcal{X}}$ be the closed subspace of $\mathcal{X}$ consisting of 
those $\mathcal{U}=(u^+,\phi^+,u^-,\phi^-,\phi)$ that satisfy
\[\mathcal{U}(0)=0,\quad
\phi^\pm=\phi\mbox{ on $\Gamma$,}\quad
L^\pm\phi^\pm=0\mbox{ in $\Omega_\pm$,}\quad
u^\pm=0\mbox{ on $\Gamma\cup\Sigma_\pm$,}\quad
\phi^\pm=0\mbox{ on $\Sigma_\pm$.}\]
Let further $\tilde{\mathcal{Y}}$ be the closed subspace of $\mathcal{Y}$ 
consisting of those $(f^+f^-,g)\in\mathcal{Y}$ that satisfy $f^-(0)=0$.
By Theorem \ref{mainthrmlin}, the linear operator 
$\mathcal{C}:\tilde{\mathcal{X}}\longrightarrow\tilde{\mathcal{Y}}$ given by
\[\mathcal{C} \mathcal{U}=(\tilde L^+u^+,\tilde L^-u^-,
\partial_t \phi - \alpha^- \phi_{z_n}^- - \alpha^+ \phi_{z_n}^+ + \zeta \cdot 
\nabla^{'} \phi - \tilde{\alpha}^- u_{z_n}^- - \tilde{\alpha}^+ u_{z_n}^+)\]
is an isomorphism. 

We rewrite system (\ref{consphi})--(\ref{demphi}) equivalently as 
\begin{equation} \label{nonlinabs} 
\mathcal{C} \mathcal{U} = \mathcal{F} (\mathcal{U}) + \mathcal{G}_0,  \qquad 
\mathcal{U} \in \tilde{\mathcal{X}},
\end{equation}

where
\begin{eqnarray*}
\mathcal{F} (\mathcal{U}) & := & (K^+ \phi^+ + R^+(\phi^+,u^+), K^- \phi^- + 
R^-(\phi^-,u^-),R^B(\mathcal{U})), \\
\mathcal{G}_0 & := & (0, 0, g_0),\\
g_0&:=&1/\mu - \partial_t\sigma+\frac{1+|\nabla'\sigma|^2}{1+\sigma_{z_n}}
(-\alpha U^-_{z_n}+\beta U^+_{z_n})
\end{eqnarray*}
Observe that $g_0 \in Y^B_\theta(\Gamma)$. Hence $\mathcal{G}_0 \in 
\tilde{\mathcal{Y}}$, and after substituting $\mathcal{V} := \mathcal{U} 
- \mathcal{F}_0$, $\mathcal{F}_0 := \mathcal{C}^{-1}\mathcal{G}_0 =: 
(0,\phi^+_0,0, \phi^-_0, \phi_0)$ our problem takes the form 
\[
\mathcal{V} = \mathcal{C}^{-1} (\mathcal{F} (\mathcal{V} +  \mathcal{F}_0)), 
\qquad \mathcal{V} \in \tilde{\mathcal{X}}.
\]

We are going to show that the mapping $\Phi: \mathcal{V} \mapsto 
\mathcal{C}^{-1} (\mathcal{F} (\mathcal{V} +  \mathcal{F}_0))$ has a fixed point 
in the closed ball $\mathbb{B}$ of radius $1$ in the space $\tilde{\mathcal{X}}$ 
provided $T > 0$ is small enough. First we make sure that $\Phi$ maps this ball 
into itself:

In view of Theorem \ref{mainthrmlin} it suffices to show that there exist $c, \delta > 0$ such that 
\[
\Vert \mathcal{F} (\mathcal{V} +  \mathcal{F}_0) \Vert_{\tilde{\mathcal{Y}}} \leq c \; T^\delta, \quad \mathcal{V} \in \mathbb{B}, \; T \leq \tau_0 \in (0,\tau), 
\]
which is implied by the estimates ($\tilde{\phi}^\pm = \phi^\pm + \phi_0^\pm$)
\begin{equation}\label{estrhs}
\begin{array}{rcl}
\Vert  K^+(\tilde{\phi}^+) + R^+(\tilde{\phi}^+,u^+) \Vert_{Y^+(\Omega_+)} & \leq & c \; T^\delta, \\
\Vert  K^-(\tilde{\phi}^-) + R^-(\tilde{\phi}^-,u^-) \Vert_{Y_\theta^-(\Omega_-)} & \leq & c \; T^\delta, \\
\Vert  R^B(\tilde{\phi}^\pm,\phi,u^\pm) \Vert_{Y_\theta^B(\Gamma)} & \leq & c \; T^\delta 
\end{array}
\end{equation}
$(\mathcal{V} \in \mathbb{B})$. Let us first  consider the parts from the elliptic phase. 

Observe that the matrices $A(p)$, $\partial_i A (p)$ have entries of the form $\frac{P(p_1,...,p_{n-1})}{(1+p_n)^j}$, where $P$ is a polynomial of degree $\leq 2$ (possibly $0$) and $j \geq 1$.

\vspace{5mm}

In view of Lemmas \ref{mult2}, \ref{mult3} we need to estimate the terms $( 1+\sigma_{z_n} + s \tilde{\phi}^-_{z_n} )^{-1}$ ($s \in [0,1]$) in the norms of $W_p^\theta((0,T),W_p^{1-\theta}(\Omega_-))$ and $L_\infty(0,T;W_p^{1-\theta}(\Omega_-))$.  

\vspace{5mm}

Recall that $1 + \sigma_{z_n}$ is invertible in the Banach algebras $H_p^1((0,\tau) \times \Omega)$, \newline  
$W_p^\theta((0,\tau),W_p^{1-\theta}(\Omega_-))$ and 
$L_\infty((0,\tau),W_p^{1-\theta}(\Omega_-))$. Since the group of invertible 
elements of a Banach algebra is open, it follows from Corollary \ref{emb} and 
$\tilde{\phi}^- \in \Mathring{X}^-(\Omega_-)$ that $f := 1 + \sigma_{z_n} + s 
\tilde{\phi}^-_{z_n}$ is invertible in $L_\infty(0,T;W_p^{1-\theta}(\Omega_-))$ 
provided $T < \tau$ is small enough.  Moreover, the inversion formula 
\[
(1 + \sigma_{z_n} + s \tilde{\phi}^-_{z_n})^{-1} =  (1+\sigma_{z_n})^{-1} 
\sum_{j=0}^\infty \big{(} s \tilde{\phi}^-_{z_n} (1+\sigma_{z_n})^{-1} 
\big{)}^j 
\]
and Corollary \ref{emb} imply that $( 1+\sigma_{z_n} + s \tilde{\phi}^-_{z_n} )^{-1}$ can be assumed to be bounded in $L_\infty(0,T;W_p^{1-\theta}(\Omega_-))$ independently of $s \in [0,1]$, $T \in (0,\tau)$ sufficiently small and \newline 
$\Vert \phi^- \Vert_{\Mathring{X}^-(\Omega_-)} \leq 1$. From 
\begin{equation} \label{uniest}
\begin{array}{cl}
& \displaystyle\int_0^T \int_0^T \frac{\Vert \frac{1}{f(t)} - \frac{1}{f(s)} 
\Vert^p_{W_p^{1-\theta}(\Omega_-)}}{|t-s|^{1+\theta p}} \; dt \, ds \\
\leq & \displaystyle\Big\| \frac{1}{f} 
\Big\|^{2p}_{L_\infty(0,T;W_p^{1-\theta}(\Omega_-))} \int_0^T \int_0^T 
\frac{\Vert f(t) - f(s) \Vert^p_{W_p^{1-\theta}(\Omega_-)}}{|t-s|^{1+\theta p}} 
\; dt \, ds
\end{array}
\end{equation}
one easily concludes that $(1 + \sigma_{z_n} + s \tilde{\phi}^-_{z_n})^{-1} \in W_p^\theta((0,T),W_p^{1-\theta}(\Omega_-))$ and that it can be assumed to be bounded in this class independently of $s \in [0,1]$, $T \in (0,\tau)$ sufficiently small and $\Vert \phi^- \Vert_{\Mathring{X}_\theta^-(\Omega_-)} \leq 1$. 

\vspace{5mm}

From this, Lemmas \ref{QV}, \ref{mult2}, \ref{mult3}, Remark \ref{regsigma} and $\phi^-, \phi_0^- \in  \Mathring{X}_{\theta'}^-$ (which is contained in $\Mathring{X}_\theta^-$) one can easily derive the estimate   
\[
\Vert  K^-(\tilde{\phi}^-) \Vert_{Y_\theta^-(\Omega_-)} \leq c \; T^\delta.
\] 
For example, we may estimate 
\begin{eqnarray*}
& & \Big\| 
\frac{U^-_{z_n}}{(1+\sigma_{z_n})^2}\tilde{\phi}^-_{z_n}\Lambda_\sigma\sigma 
\Big\|_{W_p^\theta((0,T),W_p^{-\theta}(\Omega_-))} \\
& \leq & C \; \Big[ \Big\| \frac{U^-_{z_n}}{(1+\sigma_{z_n})^2} 
\Big\|_{L_\infty(0,T;W_p^{1-\theta}(\Omega_-))} \Vert 
\tilde{\phi}^-_{z_n}\Lambda_\sigma\sigma 
\Vert_{W_p^\theta((0,T),W_p^{-\theta}(\Omega_-))} \\
& & + \Big\| \frac{U^-_{z_n}}{(1+\sigma_{z_n})^2}  
\Big\|_{W_p^\theta((0,T),W_p^{1-\theta}(\Omega_-))} \Vert 
\tilde{\phi}^-_{z_n}\Lambda_\sigma\sigma 
\Vert_{L_\infty(0,T;W_p^{-\theta}(\Omega_-))} \Big]
\end{eqnarray*}
(Lemmas \ref{mult2}, \ref{mult3} Remark \ref{regsigma}),  
\begin{eqnarray*}
\Vert \tilde{\phi}^-_{z_n}\Lambda_\sigma\sigma \Vert_{L_\infty(0,T;W_p^{-\theta}(\Omega_-))} & \leq & C \; \Vert \tilde{\phi}^-_{z_n}\Lambda_\sigma\sigma \Vert_{W_p^\theta((0,T),W_p^{-\theta}(\Omega_-))} \\
& \leq & C \; \big{[} \Vert \Lambda_\sigma\sigma \Vert_{L_\infty(0,\tau;W_p^{-\theta}(\Omega_-))} \Vert \tilde{\phi}^-_{z_n} \Vert_{W_p^\theta((0,T),W_p^{1-\theta'}(\Omega_-))} \\
& & + \Vert \tilde{\phi}^-_{z_n} \Vert_{L_\infty(0,T;W_p^{1-\theta}(\Omega_-))} \Vert \Lambda_\sigma\sigma \Vert_{W_p^\theta((0,\tau),W_p^{-\theta}(\Omega_-))} \big{]} 
\end{eqnarray*}
(Lemma \ref{mult2}, \ref{mult3}, Remark \ref{regsigma}) and 
\begin{eqnarray*}
\Vert \tilde{\phi}^-_{z_n} \Vert_{L_\infty(0,T;W_p^{1-\theta}(\Omega_-))} & \leq & c \; T^\delta, \qquad \delta \in (0,\theta - 1/p), \\
\Vert \tilde{\phi}^-_{z_n} \Vert_{W_p^\theta((0,T),W_p^{1-\theta'}(\Omega_-))} & \leq & c \; (T^{1/p+\delta} + T^{\theta'-\theta})
\end{eqnarray*}
(Corollary \ref{emb}), since $\tilde{\phi}^-(0) = 0$.

\vspace{3mm}

By construction, $R^-(\tilde{\phi}^-,u^-)$ and $R^B(\tilde{\phi}_0^\pm,\phi,u^\pm)$ are quadratic in terms that vanish at $t=0$. Hence, the second and third inequality in (\ref{estrhs}) are easy consequences of Lemmas \ref{mult2}, \ref{mult3} and \ref{alg} (cf. also Remarks \ref{rem1}, \ref{rem2}). 

\vspace{5mm}

The 'parabolic parts' of (\ref{estrhs}) can be treated along the lines of \cite{solfr}. We restrict ourselves to give a remark on the invertibility of the terms $1 + \sigma_{z_n} + s \tilde{\phi}^+_{z_n}$ in $L_\infty((0,T) \times \Omega_+)$: Since $p > n + \frac{3}{2} + \sqrt{2n + \frac{1}{4}} > n + 2$ there is a $\sigma \in (\frac{1}{p}, \frac{1 - N/p}{2})$. We choose $\delta > 0$ such that $\sigma > \frac{1}{p} + \delta$. Then
\[
\tilde{\phi}^+_{z_n} \in \Mathring{W}_p^\sigma((0,T),W_p^{1-2\sigma}(\Omega_+)) \hookrightarrow \Mathring{W}_p^\sigma((0,T),C(\bar{\Omega}_+)).
\]
Hence the 'uniform invertibility' of $1 + \sigma_{z_n} + s \tilde{\phi}^+_{z_n}$ follows from Corollaries \ref{emb}, \ref{embpara}. 

\vspace{5mm}

The same kind of arguments (oriented again at \cite{solfr}) show that there exist $c, \delta > 0$ such that 
\[
\Vert \Phi(\mathcal{V}_1) - \Phi(\mathcal{V}_2) \Vert_{\mathcal{X}} \leq c \; T^\delta \; \Vert \mathcal{V}_1 - \mathcal{V}_2 \Vert_{\mathcal{X}} \quad (\mathcal{V}_1, \mathcal{V}_2 \in \mathbb{B}),
\]
showing that $\Phi$ is a contraction on $\mathbb{B}$ provided $T > 0$ is small enough. This finally proves Theorem \ref{mainthrm}.

\appendix
\section{Appendix}
\subsection{Parameter dependent elliptic problems}

\begin{remark}\label{absrem}
{\rm Let $X_i$, $Y_i$, $i=0,1$ be Banach spaces, $q\in(1,\infty)$, $J=(0,T)$. 
Assume $Y_1\hookrightarrow Y_0$, 
\[\begin{array}{rl}
A\in& L_\infty(J,\cll_{is}(X_1,Y_1)), \\
A^{-1}\in& L_\infty(J,\cll_{is}(Y_1,X_1)),\\
K\in& L_q(J,\cll(X_0,Y_1))\;\cap\; L_\infty(J,\cll(X_0,Y_0)),\\
A+K\in& L_\infty(J,\cll_{is}(X_0,Y_0)),\\
(A+K)^{-1}\in& L_\infty(J,\cll_{is}(Y_0,X_0)),\\
F\in& L_\infty(J,Y_0)\;\cap\; L_q(J,Y_1).
\end{array}\]
and define $u\in L_\infty(J,X_0)$ by 
\[u(t):=(A(t)+K(t))^{-1}F(t).\]
Then $u(t)=A(t)^{-1}(F(t)-K(t)u(t))$ and therefore $u\in L_q(J,X_1)$, with estimate
\bea
\|u\|_{L_q(J,X_1)}&\leq& 
\|A^{-1}\|_{L_\infty(J,\cll_{is}(Y_1,X_1))}(\|F\|_{L_q(J,Y_1)}
+\|K\|_{L_q(J,\cll(X_0,Y_1))}\|u\|_{L_\infty(J,X_0)})\\
&\leq&\|A^{-1}\|_{L_\infty(J,\cll_{is}(X_1,Y_1))}(\|F\|_{L_q(J,Y_1)}\\
&&+\|K\|_{L_q(J,\cll(X_0,Y_1))} \|(A+K)^{-1}\|_{L_\infty(J,\cll_{is}(Y_0,X_0))}
\|F\|_{L_\infty(J,Y_0)}).
\eea }
\end{remark}

\begin{lemma}\label{hoelder}
Let $0<\beta<\alpha<1$,  $\subset\RRM^n$ be a 
domain with 
$d:=\diam D \leq 1$ and $\phi\in BUC^\alpha(D)$ such that $\phi(x_0)=0$ 
for some $x_0 \in D$. Then
\[\|\phi\|_{BUC^{\beta}(D)}\leq 
d^{\alpha-\beta}\|\phi\|_{BUC^{\alpha}(D)}.\]
\end{lemma}
{\bf Proof:}
We have for $x,y \in D$
\bea
|\phi(x)|&=&|\phi(x)-\phi(x_0)|\leq d^\alpha\|\phi\|_{BUC^{\alpha}(D)},\\
|\phi(x)-\phi(y)|&\leq&|x-y|^\alpha\|\phi\|_{BUC^{\alpha}(D)}
\leq d^{\alpha-\beta}|x-y|^\beta\|\phi\|_{BUC^{\alpha}(D)}.
\eea
\qed

\begin{lemma}\label{remain}
(H\"older estimates for coefficients in remainder terms)

Let 
\[0<\theta_1<\theta_2<\theta_3
<\frac{1}{2}\left(1-\frac{n+1}{p}\right),\]
$a\in Z^0(\Omega)$, $\xi_0\in\Omega$, $a(\xi_0,0)=0$, $\lambda\in (0,1)$,
$\tau_0 \in J$. Let 
$\Omega_\lambda:=\Omega_-\cap B(\xi_0,\lambda) \neq \emptyset$, 
$Q_{\lambda,\tau_0}:=\Omega_\lambda\times(0,\tau_0)$.
Then
\[a|_{Q_{\lambda,\tau_0}}\in
BUC^{\theta_1}((0,\tau_0),BUC^{\theta_2}(\Omega_\lambda))\]
and
\[\|a|_{Q_{\lambda,\tau_0}}\|_{BUC^{\theta_1}((0,\tau_0),BUC^{\theta_2}
(\Omega_\lambda))}
\leq 
C(\tau_0^{\theta_2-\theta_1}+\lambda^{2\theta_3-\theta_2})
\|a\|_{Z^0(\Omega)}.\]
\end{lemma}
{\bf Proof:}  We rewrite
\[a(x,t)=\hat a(x,t)+a(x,0),\qquad x\in\Omega,\,t\in J\]
and estimate the terms on the right separately. (The second term will be
interpreted both as a function on $\Omega$ and as a function on $\Omega\times
J$ which is constant with respect to $t$.)
Trace and embedding theorems yield $a(\cdot,0)\in BUC^{2\theta_3}(\Omega)$, 
\[\|a(\cdot,0)\|_{BUC^{2\theta_3}(\Omega)}\leq C\|a\|_{Z^0(\Omega)},\]
and thus by Lemma \ref{hoelder}
\be\label{esta1}
\|a(\cdot,0)|_{\Omega_\lambda}\|_{BUC^{\theta_2}(\Omega_\lambda)}
\leq C\lambda^{2\theta_3-\theta_2}\|a\|_{Z^0(\Omega)}.
\ee
Furthermore, as $2\theta_2<1-\frac{n+1}{p}$ we have by embedding
\[\hat a\in BUC^{\theta_2}(J,BUC^{\theta_2}(\Omega))\]
and by restriction
\[\hat a|_{Q_{\lambda,\tau_0}}\in
BUC^{\theta_2}((0,\tau_0),BUC^{\theta_2}(\Omega_\lambda))\]
with estimate
\[\|\hat a|_{Q_{\lambda,\tau_0}}
\|_{BUC^{\theta_2}((0,\tau_0),BUC^{\theta_2}(\Omega_\lambda))}
\leq C\|a\|_{Z^0(\Omega)}.\]
Consequently, using $\hat a|_{t=0}=0$,
\[\|\hat a|_{Q_{\lambda,\tau_0}}
\|_{BUC^{\theta_1}((0,\tau_0),BUC^{\theta_2}(\Omega_\lambda))}
\leq C\tau_0^{\theta_2-\theta_1}\|a\|_{Z^0(\Omega)}.\]
Together with (\ref{esta1}) this implies the result. \qed

\bigskip

Let $a_{ij}\in Z^0(\Omega)$ be uniformly elliptic, i.e. 
\be\label{ellcoeff}
a_{ij}\xi_i\xi_j\geq\mu|\xi|^2 \mbox{ in $\Omega \times J$}
\ee
for some $\mu>0$. Let $b_i\in Y^-(\Omega_-)$, $f\in Y^-(\Omega_-)$, $g\in X^-_{\tr}(\Omega_-)$.
\begin{lemma}\label{elleq} 
Let $M=\max\{\|a_{ij}\|_{Z^0},\|b_i\|_{Y^-},\mu^{-1}\}.$ There are 
constants $C(M)$, $T_0(M)$ such that for $T\leq T_0$, there is precisely one 
solution $u\in X^-$ to the time-dependent elliptic 
problem
\[\left.
\begin{array}{rcll}
a_{ij}\partial_{ij}u+b_i\partial_i u&=&f&\mbox{ in $\Omega_-\times J$,}\\
u&=&g&\mbox{ on $\Gamma\times J$,}\\
u&=&0&\mbox{ on $\Sigma_-\times J$.}
\end{array}\right\}\]
It satisfies
\[\|u\|_{X^-(\Omega_-)}\leq C(M)(\|f\|_{Y^-(\Omega_-)}+\|g\|_{X^-_\tr(\Omega_-)})\]
If $g=0$ and $f\in\Mathring{Y}^-(\Omega_-)$ then $u\in\Mathring{X}^-(\Omega_-)$. 
\end{lemma}
{\bf Proof:}
1. We first show $u\in L_p(J,H^2_p(\Omega_-))$ with the corresponding estimate. For this 
we set 
\bea
X_0&:=&W_p^{2-\theta}(\Omega_-),\\
X_1&:=&H^2_p(\Omega_-),\\
Y_0&:=&W_p^{-\theta}(\Omega_-)\times W_p^{2-\theta-1/p}(\Gamma)
\times W_p^{2-\theta-1/p}(\Sigma_-),\\
Y_1&:=&L_p(\Omega_-)\times W_p^{2-1/p}(\Gamma)
\times W_p^{2-1/p}(\Sigma_-),\\
A&:=&(a_{ij}\partial_{ij},\Tr_\Gamma,\Tr_\Sigma),\\
K&:=&(b_i\partial_i,0,0),\\
F&:=&(f,g,0)
\eea
and aim at the application of Remark \ref{absrem}. $F$ clearly satisfies the 
assumptions.

1.1. For $v\in X_0$, $t\in J$ we have
\[\|K(t)v\|_{Y_0}=\|b_i\partial_iv\|_{W^{-\theta}_p(\Omega_-)}
\leq C\|\vec b\|_{Y^-}\|v\|_{X_0}\]
due to the embedding $X_0\hookrightarrow BUC^{1+\theta'}(\Omega_-)$ for a 
suitable $\theta'\in(\theta,1)$ and the pointwise multiplicator property
(\cite{trfs1} Theorem 3.3.2)
\be\label{mult}
BUC^{\theta'}(\Omega_-)\cdot W_p^{-\theta}(\Omega_-)\hookrightarrow
W_p^{-\theta}(\Omega_-). 
\ee
So $K\in L_\infty(J,\cll(X_0,Y_0))$.
 Furthermore,
\bea 
&&\int_J\|K(t)\|_{\cll(X_0,Y_1)}^p\,dt=\int_J\left(
\sup_{\|v\|_{X_0}=1}\|b_i\partial_i v\|_{L_p(\Omega_-)}
\right)^p\,dt\\
&=&\int_J\sup_{\|v\|_{X_0}=1}\int_{\Omega_-}|b_i\partial_i v|^p\,dxdt
\leq C\|\vec b\|^p_{Y^-(\Omega_-)},
\eea
where we use the embedding $X_0\hookrightarrow BUC^1(\Omega_-)$.
This shows $K\in L_p(\cll(X_0,Y_1))$.

1.2. By parallel reasonings, we get $A\in BUC(J,\cll(X_1,Y_1))$,
$A+K\in BUC(J,\cll(X_0,Y_0))$ with norms depending only on 
$\|a_{ij}\|_{Z^0(\Omega)}$, $\|b_i\|_{Y^-(\Omega_-)}$. The fact that for $t\in J$ we have
$A(t)\in\cll_{is}(X_1,Y_1)$, with $\|A(t)^{-1}\|_{\cll(Y_1,X_1)}$ depending 
only on $\|a_{ij}\|_{Z^0(\Omega)}$ and $\mu$ follows from standard theory on 
elliptic boundary value problems. To get $A(t)+K(t)\in \cll_{is}(X_0,Y_0)$ one 
proceeds as in the proof of \cite{trfs1} Theorem 4.3.3., with slight 
modifications due to the fact that the coefficients of $A$ and $K$ are not 
$C^\infty$. The proof remains valid, anyway, as by (\ref{mult}) and 
interpolation we have estimates for the lower order term of the type
\begin{eqnarray*}
\|b_i(t)\partial_iw\|_{W_p^{-\theta}(\Omega_-)} & \leq & C\|\vec b\|_{Y^-(\Omega_-)}
\|w\|_{BUC^{1+\theta'}(\Omega_-)} \\
& \leq & \eps\|w\|_{X_0}+C(\eps,\|\vec b\|_{Y^-(\Omega_-)})
\|w\|_{W_p^{-\theta}(\Omega_-)}
\end{eqnarray*}
for any $\eps>0$. Note also that we have to use (\ref{mult}) together  with 
Lemma \ref{hoelder} (with $\beta=\theta$, $\alpha=\theta'$) to estimate the 
highest-order error terms occurring from freezing of the coefficients. 

As all assumptions of the above remark are valid, we conclude 
$u\in L_p(J,H^2_p(\Omega_-))$ and
\[\|u\|_{L_p(J,H^2_p(\Omega_-))}\leq C(\|f\|_{Y^-(\Omega_-)}+\|g\|_{X^-_\tr(\Omega_-)}).\]

2. By arguments as above, we have
\bea
&&\int_J\int_J\frac{\|(K(t)-K(s))\|_{\cll(X_0,Y_0)}^p}{|t-s|^{1+p\theta}}\;
ds\,dt\\
&=&\int_J\int_J\frac{\Big(\displaystyle\sup_{\|v\|_{X_0}=1}
\|(b_i(t)-b_i(s))\partial_iv\|_{W^{-\theta}_p(\Omega_-)}\Big)^p}{|t-s|^{ 
1+p\theta}}\;ds\,dt\\
&\leq& C\int_J\int_J\frac{\|\vec b(t)-\vec 
b(s)\|_{W^{-\theta}_p(\Omega_-)}^p}{|t-s|^{1+p\theta}}\sup_{\|v\|_{X_0}=1}
\|v\|_{BUC^{1+\theta'}(\Omega_-)}^p\;ds\,dt\leq C\|\vec b\|_{Y^-(\Omega_-)}^p,
\eea
hence $K\in W_p^{\theta}(J,\cll(X_0,Y_0))$ and by parallel arguments $A\in
W_p^{\theta}(J,\cll(X_0,Y_0))$, with norms depending on $\|a_{ij}\|_{Z^0(\Omega)}$, 
$\|b_i\|_{Y^-(\Omega_-)}$. Write $\hat A:=A+K$. By Step 1.2, $\hat A(0)\in 
\cll_{is}(X_0.Y_0)$, and therefore there is a $\delta>0$ such that the open 
ball $\clb$ in $\cll(X_0,Y_0)$ centered at $\hat A(0)$ with radius $\delta$ lies 
in $\cll_{is}(X_0,Y_0)$, and the operator $\inv$ given by $\inv(B)=B^{-1}$ is 
in $BUC^1(\clb,\cll(Y_0,X_0))$. Because $\hat A$ is continuous in time 
with values in $\cll(X_0,Y_0)$, by shrinking $J$ if necessary, we can arrange 
that $\hat A(t)\in\clb$ for all $t\in J$. Therefore we have 
\[\|\hat A(t)^{-1}-\hat A(s)^{-1}\|_{\cll(Y_0,X_0)}
\leq C\|\hat A(t)-\hat A(s)\|_{\cll(X_0,Y_0)},\]
and from this it follows straightforwardly that 
$\inv\circ\hat A\in W_p^\theta(J,\cll(Y_0,X_0))$. Finally from 
this  and 
$F\in W_p^\theta(J,Y_0)$ we get $u\in 
W_p^\theta(J,W_p^{-\theta}(\Omega_-))$ and
\[\|u\|_{W_p^\theta(J,W_p^{-\theta}(\Omega_-))}
\leq C(\|f\|_{Y^-(\Omega_-)}+\|g\|_{X^-_\tr(\Omega_-)}).\]
 \qed

\subsection{Uniform anisotropic embeddings, multiplications and related estimates}
\begin{lemma} \label{fort}
Let $m \in \mathbb{N}$, $q \in (1,\infty)$, $\sigma_1 = 0$ or $\sigma_1 \in (1/q, 1]$. Let further $1/q < \sigma_2 \leq \cdots \leq \sigma_m \leq 1$ and  $\sigma_1 \leq \sigma_2$. Let $X_1,...,X_m$ be Banach spaces such that $X_1 \hookrightarrow X_2 \hookrightarrow \cdots \hookrightarrow X_m$. Given $T \in (0,\infty]$, let 
\[
W(T) := \bigcap_{i = 1}^m W_q^{\sigma_i}((0,T),X_i), \qquad W:= W(\infty).
\]
There exist bounded linear operators $\mathcal{E}_T: W(T) \rightarrow W$ and a constant $C > 0$ such that
\[
\Vert \mathcal{E}_T u \Vert_{L_q(0,\infty;X_i)} \leq C \Vert u \Vert_{L_q(0,T;X_i)}
\]
for all $T \in (0,\infty]$, $u \in L_q(0,T;X_i)$ and $i = 1,...,m$, and (in case that either $\sigma_1 \in (1/q, 1]$ or $\sigma_1 = 0$ and $m \geq 2$)
\[
\Vert \mathcal{E}_T u \Vert_W \leq C \Vert u \Vert_{W(T)} 
\]
for all $T \in (0,\infty]$ and $u \in \Mathring{W}(T)$.
\end{lemma}

{\bf Proof:}
For $m = 1$ and $0 < T \leq T_0$ this is stated and proved in Proposition $6.1$ in \cite{pss}. The case of a general $m \in \mathbb{N}$ is an immediate consequence. Moreover, a careful inspection of the proof shows that the constant $C$ can be chosen independent of $T > T_0$, too.
\qed

\begin{corollary} \label{emb}
Let $q \in (1,\infty)$, $\delta \geq 0$, $1/q + \delta < \sigma \leq 1$ and let $X$ be a Banach space. Given $T \in (0,\infty]$ we have $W_q^\sigma((0,T),X) \hookrightarrow BUC^\delta((0,T),X)$. There exists a constant $C > 0$ such that 
\[
\Vert u \Vert_{BUC^\delta((0,T),X)} \leq C \Vert u \Vert_{W_q^\sigma((0,T),X)}
\]
for all $T \in (0,\infty]$ and $u \in \Mathring{W}_q^\sigma((0,T),X)$.
\end{corollary}

\begin{corollary} \label{embpara}
Let $q \in (1,\infty)$, $2 \delta + \sigma \leq 2$. Given $T \in (0,\infty]$ we have 
\[
W_q^1((0,T),W_q^2(\Omega_+)) \cap L_q((0,T),W_q^2(\Omega_+)) \hookrightarrow W_q^\delta((0,T),W_q^\sigma(\Omega_+)). 
\]
There exists a constant $C > 0$ such that 
\[
\Vert u \Vert_{W_q^\delta((0,T),W_q^\sigma(\Omega_+))} \leq C \Vert u \Vert_{W_q^1((0,T),W_q^2(\Omega_+)) \cap L_q((0,T),W_q^2(\Omega_+))}
\]
for all $T \in (0,\infty]$ and $u \in \Mathring{W}_q^1((0,T),W_q^2(\Omega_+)) \cap L_q((0,T),W_q^2(\Omega_+))$.
\end{corollary}

\vspace{5mm}

Let $X,Y,Z$ be Banach spaces whose elements can be interpreted as real-valued
functions on the same domain of definition. Then there is a pointwise product
$(u,v)\mapsto uv$ on $X\times Y$. We write
\[X \cdot Y \hookrightarrow Z\]
if for all $(u,v)\in X \times Y$ we have $uv \in Z$ and there is an $M>0$ such that
\[\|uv\|_Z \leq M\|u\|_X\|v\|_Y,\qquad(u,v)\in X\times Y.\]

\vspace{5mm}

\begin{lemma} \label{mult2}
Let $q \in (1,\infty)$, $1 > \rho > \sigma > 1/q$, $T > 0$ and let $X,Y,Z$ be Banach spaces s.t. $X \cdot Y \hookrightarrow Z$. The following holds true:
\begin{itemize}
\item[i)] $C([0,T],X) \cdot L_q(0,T;Y) \hookrightarrow L_q(0,T;Z)$ and 
\[
\Vert u v \Vert_{L_q(0,T;Z)} \leq M \Vert u \Vert_{C([0,T],X)} \Vert v \Vert_{L_q(0,T;Y)}
\]
for all $u \in C([0,T],X)$, $v \in L_q(0,T;Y)$. 
\item[ii)] $C^\rho([0,T],X) \cdot W_q^\sigma((0,T),Y) \hookrightarrow W_q^\sigma((0,T),Z)$ and 
\begin{eqnarray*}
\Vert uv \Vert_{W_q^\sigma((0,T),Z)} & \leq & C(\rho,\sigma,q,M) \; \big[ \Vert u \Vert_{L_\infty(0,T;X)} \Vert v \Vert_{W_q^\sigma((0,T),Y)} \\ 
& & + \, T^{\rho-\sigma + 1/q} \Vert v \Vert_{L_\infty(0,T;Y)} \Vert u \Vert_{C^{\rho}([0,T],X)} \big]
\end{eqnarray*}
for all $u \in C^{\rho}([0,T],X)$, $v \in W_q^\sigma((0,T),Y)$. 
\item[iii)] $W_q^\sigma((0,T),X) \cdot W_q^\sigma((0,T),Y) \hookrightarrow W_q^\sigma((0,T),Z)$ and 
\begin{eqnarray*}
\Vert uv \Vert_{W_q^\sigma((0,T),Z)} & \leq & C(q,M) \; \big[ \Vert u \Vert_{L_\infty(0,T;X)} \Vert v \Vert_{W_q^\sigma((0,T),Y)} \\ 
& & + \, \Vert v \Vert_{L_\infty(0,T;Y)} \Vert u \Vert_{W_q^\sigma((0,T),X)} \big]
\end{eqnarray*}
for all $u \in W_q^\sigma((0,T),X)$, $v \in W_q^\sigma((0,T),Y)$. 
\end{itemize}
\end{lemma}

{\bf Proof:}
The first statement is trivial. Let $u \in C^\rho([0,T],X)$, $v \in W_q^\sigma((0,T),Y)$. We have  
\[
\int_0^T \Vert u(t) v(t) \Vert_Z^q \; dt \leq M \Vert u \Vert_{L_\infty(0,T;X)}^q \int_0^T \Vert v(t) \Vert_Y^q \; dt
\]
and
\begin{eqnarray} \label{est1}
\int_0^T \int_0^T \frac{\Vert u(t)v(t) - u(s)v(s) \Vert_Z^q}{|t-s|^{1+\sigma q}} \; dt \, ds & \leq & (2M)^q \; \big[ \; \int_0^T \int_0^T \frac{\Vert u(t) \Vert_X^q \Vert v(t)-v(s) \Vert_Y^q}{|t-s|^{1+\sigma q}} \; dt \, ds \nonumber \\
& & + \; \int_0^T \int_0^T \frac{\Vert v(s) \Vert_Y^q \Vert u(t)-u(s) \Vert_X^q}{|t-s|^{1+\sigma q}} \; d(t,s) \; \big] \nonumber \\
& \leq & (2M)^q \; \big[ \; \Vert u \Vert_{L_\infty(0,T;X)}^q \Vert v \Vert_{W_q^\sigma((0,T),Y)}^q \nonumber \\
& & + \; \Vert v \Vert_{L_\infty(0,T;Y)}^q \Vert u \Vert_{C^\rho([0,T],X)}^q \\
& & \times \; \int_0^T \int_0^T |t-s|^{(\rho-\sigma) q - 1} \; dt \, ds \; \big]. \nonumber 
\end{eqnarray}
All our assertions follow easily from these estimates.
\qed

\vspace{5mm}

An immediate consequence is the following  

\begin{lemma} \label{mult3}
Under the assumptions of Lemma {\rm \ref{mult2}} we have
\begin{eqnarray*}
\Vert uv \Vert_{W_q^\sigma((0,T),Z)} & \leq & C(\rho,\sigma,q,M) \; \big[ T^\rho \; \Vert u \Vert_{C^{\rho}([0,T],X)} \Vert v \Vert_{W_q^\sigma((0,T),Y)}  \\ 
& & + \, T^{\rho+\varepsilon-\sigma+1/q} \; \Vert v \Vert_{W_q^\sigma((0,T),Y)} \Vert u \Vert_{C^{\rho}([0,T],X)} \big]
\end{eqnarray*}
for all $u \in \Mathring{C}^{\rho}([0,T],X)$, $v \in \Mathring{W}_q^\sigma((0,T),Y)$ and $0 \leq \varepsilon < \sigma - 1/q$. Moreover, 
\[
\Vert uv \Vert_{W_q^\sigma((0,T),Z)} \leq C(q,M) \; T^\varepsilon \; \Vert u \Vert_{W_q^\sigma((0,T),X)} \Vert v \Vert_{W_q^\sigma((0,T),Y)}
\]
for all $u \in \Mathring{W}_q^\sigma((0,T),X)$, $v \in \Mathring{W}_q^\sigma((0,T),Y)$ and $0 \leq \varepsilon < \sigma-1/q$.
\end{lemma}

\begin{remark} \label{rem1} {\rm [product estimate, elliptic phase]
Let $q \in (1,\infty)$, $N \in \mathbb{N}$ and let $D \subset \RRM^N$ be open. Lemmas {\rm \ref{mult2}} and {\rm \ref{mult3}} guarantee smallness of terms 
\[
\Vert D^2 u D v \Vert_{L_q(0,T;L_q(D)) \cap W_q^\sigma((0,T),W_q^{-\sigma}(D))} 
\]
for small values of $T$ and $u,v \in L_q(0,T;W_q^2(D)) \cap \Mathring{W}_q^\sigma((0,T),W_q^{2-\sigma}(D))$ by choosing  $Y=Z=W_q^{-\sigma}(D)$, $X = W_q^{1-\sigma}(D)$ ($1/q < \sigma < 1/2$) and $Y=Z=L_q(D)$, $X = W_q^{1-\sigma}(D)$ ($1/q < \sigma < 1-N/q$), respectively. 

\vspace{5mm}

Observe that the conditions $1/q < \sigma < 1/2$ and $1/q < \sigma < 1-N/q$ are both satisfied if $\frac{1}{q}\frac{q-1}{q-N} < \sigma < \frac{1}{2}(1-\frac{N+1}{q})$ and $q > N-1$.}
\end{remark}

\begin{lemma} \label{alg}
Let $N \in \mathbb{N}$, $1 > \sigma > 1/q$, $1 \geq r > \frac{N-1}{q-1/\sigma}$ and let $U$ be an open set in $\RRM^{N-1}$. For $T \in (0,\infty]$ let  
\[ E(T) := L_q(0,T;W_q^r(U)) \cap W_q^\sigma((0,T),L_q(U)). \]  
Then $E(T) \hookrightarrow BUC((0,T) \times U)$ and $E(T)$ is a Banach algebra. There exists a constant $C > 0$ such that   
\begin{itemize}
\item[i)]  $\Vert u \Vert_\infty \leq C \Vert u \Vert_{E(T)}$ for all $T \in (0,\infty]$ and $u \in \Mathring{E}(T)$;
\item[ii)] $\Vert uv \Vert_{E(T)} \leq C ( \Vert u \Vert_\infty \Vert v \Vert_{E(T)} + \Vert v \Vert_\infty \Vert u \Vert_{E(T)} )$ for all $T \in (0,\infty]$ and $u,v \in \Mathring{E}(T)$
\end{itemize}
(where $\Vert \cdot \Vert_\infty$ denotes the sup of a function over the set $(0,T) \times M$). 

\vspace{5mm}

If $\delta \in [0,1]$, $r > \frac{N-1}{q(1-\delta)}$ and $\delta \sigma > 1/q + \varepsilon$, then $E(T) \hookrightarrow BUC^{\varepsilon}((0,T),BUC(U))$ and there is a constant $C > 0$ such that
\begin{itemize}
\item[iii)] $\Vert u \Vert_\infty \leq C T^{\varepsilon} \Vert u \Vert_{E(T)}$
\end{itemize}
for all $T \in (0,\infty]$ and $u \in \Mathring{E}(T)$.
\end{lemma}

{\bf Proof:}
The embedding $E(T) \hookrightarrow BUC((0,T) \times U)$ is stated and proved in Lemma $4.4$ in \cite{dss} and the estimate i) follows straightforwardly from Lemma \ref{fort}, Corollary \ref{emb}. Observe that for a.e. $t \in (0,T)$ 
\[
\Vert u(t) v(t) \Vert_{W_q^r(D)} \leq C(q) \big( \; \Vert u(t) \Vert_{C(\bar{U})} \Vert v(t) \Vert_{W_q^r(U)} + \Vert v(t) \Vert_{C(\bar{U})} \Vert u(t) \Vert_{W_q^r(U)} \; \big), 
\]
since $r > \frac{N-1}{q-1/\sigma} > \frac{N-1}{q}$. Hence, 
\[
\begin{array}{rcl}
\int_0^T \Vert u(t)v(t) \Vert_{W_q^r(U)}^q \; dt & \leq & C(q) \; \big( \; \Vert u \Vert_\infty^q \Vert v \Vert_{L_q(0,T;W_q^r(U))}^q \\ 
& & + \Vert v \Vert_\infty^q \Vert u \Vert_{L_q(0,T;W_q^r(U))}^q \; \big)
\end{array}
\]
and, as calculations similar to (\ref{est1}) show, 
\[
[uv]_{\sigma,q,L_q(U)}^q \leq C(q) \big( \; \Vert u \Vert_\infty^q \Vert v \Vert_{W_q^\sigma((0,T),L_q(U))}^q + \Vert v \Vert_\infty^q \Vert u \Vert_{W_q^\sigma((0,T),L_q(U))}^q \; \big)
\]
Assertion ii) is now an easy consequence of this and of Lemma \ref{fort}, Corollary \ref{emb}. Assertion iii) follows from Lemma $4.3$ in \cite{dss} and again Lemma \ref{fort}, Corollary \ref{emb}.
\qed

\begin{remark} \label{rem2}
{\rm For $r=1-1/q$ the conditions $r > \frac{N-1}{q-1/\sigma}$, $\sigma > 1/q$ are satisfied if $\sigma > \frac{1}{q} \frac{q-1}{q-N}$. In this case, $1/(q\sigma) < (q-N)/(q-1)$. If $\delta \in (1/(q\sigma), (q-N)/(q-1))$, we have $1-1/q > (N-1)/(q(1-\delta))$ and $\delta \sigma > 1/q$. Thus, (identifying $\Gamma$ with $\RRM^{n-1}$) Lemma {\rm \ref{alg}} applies to the space $Y_\theta^B(\Gamma)$ frequently used in this paper.}
\end{remark}
 
\subsection{Some auxiliary results concerning localizations} \label{socf}
Let $r \in (0,1)$, $q \in [1,\infty)$, and let $\{ \Omega^{(k)} \}_{k \in 
\mathcal{K}}$ be the collection of sets defined in the proof of Lemma 
\ref{varcoeff}. Suppose further that 
\begin{itemize}
\item $f \in C^\infty(\Gamma)$, $\{ f_k \}, \{ g_k \} \subset C^\infty(\Gamma)$, 
$\mbox{supp}(f_k) \subset \Omega^{(k)}$ ($k \in \mathcal{K}$);
\item $\{ \psi_k \} \subset C^\infty(\Gamma)$ are such that $\mbox{supp}(\psi_k) 
\subset \Omega^{(k)}$ and $\vert \partial^\alpha \psi_k \vert_\infty \leq C 
\lambda^{-\alpha}$ uniformly for $k \in \mathcal{K}$.
%\item $\{ u_k \} \subset W_p^{1+r}(\Gamma)$ ($k \in \mathcal{K}$).
\end{itemize} 
\begin{lemma} \label{suml1}
We have 
\begin{equation} \label{sum1}
\Big\| \sum_{k \in \mathcal{K}} f_k \Big\|_{L_q(\Gamma)}^q \leq (N_0)^q \sum_{k 
\in \mathcal{K}} \Vert f_k \Vert_{L_q(\Gamma)}^q 
\end{equation}
and 
\begin{equation} \label{sum3}
\Big\| \sum_{k \in \mathcal{K}} f_k \Big\|_{W^r_q(\Gamma)}^q \leq 2(2 N_0)^q 
\sum_{k \in \mathcal{K}} \Vert f_k \Vert_{W^r_q(\Gamma)}^q. 
\end{equation}
\end{lemma}
{\bf Proof:}

1. Let $x \in \Gamma$. Since $x$ is an element of at most $N_0$ of the sets 
$\Omega^{(k)}$, the sum $\sum_{k \in \mathcal{K}} |f_k(x)|$ has at most $N_0$ 
nonzero summands. Hence 
\[
\Big{(} \sum_{k \in \mathcal{K}} |f_k(x)| \Big{)}^q \leq (N_0)^q \sum_{k \in 
\mathcal{K}} |f_k(x)|^q.
\]

2. Let $(x,y) \in \Gamma \times \Gamma$. Then the sum $\sum_{k \in \mathcal{K}} |f_k(x) - f_k(y)|$ has at most $2 N_0$ nonzero summands. Hence 
\[
\Big{(} \sum_{k \in \mathcal{K}} |f_k(x) - f_k(y)| \Big{)}^q \leq (2N_0)^q 
\sum_{k \in \mathcal{K}} |f_k(x) - f_k(y)|^q.
\]
The assertion follows from the definition of the intrinsic norms 
\begin{eqnarray*}
\Vert f_k \Vert_{W^r_q(\Gamma)} & := & \Vert f_k \Vert_{L_q(\Gamma)} + [f_k]^q_{q,r,\Gamma} \\ 
& := & \Vert f_k \Vert_{L_q(\Gamma)} + \int_{\Gamma} \int_{\Gamma} \frac{|f_k 
(x)-f_k (y)|^q}{|x-y|^{n-1+rq}} \; d\sigma(x)\;d\sigma(y).
\end{eqnarray*}
\qed
\begin{remark} \label{suml2}
{\rm A special case of Lemma \ref{suml1} are the estimates  
\begin{equation} \label{sum2}
\Big\| \sum_{k \in \mathcal{K}} \psi_k f \Big\|_{L_q(\Gamma)}^q \leq (N_0)^q 
\sum_{k \in \mathcal{K}} \Vert \psi_k f \Vert_{L_q(\Gamma)}^q 
\end{equation}
and
\begin{equation} \label{sum4}
\Big\| \sum_{k \in \mathcal{K}} \psi_k f \Big\|_{W^r_q(\Gamma)}^q \leq 2(2 
N_0)^q \sum_{k \in \mathcal{K}} \Vert \psi_k f  \Vert_{W^r_q(\Gamma)}^q. 
\end{equation} }
\end{remark}
A direct consequence of Lemma \ref{suml1} and a standard approximation argument is  
\begin{corollary} \label{suml3}
Let $V_\theta^B(\Gamma)$, $V \in \{ X,Y \}$ and $p$ be as in Section {\rm 4.1}. Then 
\begin{equation} \label{wt}
\Big\| \sum_{k \in \mathcal{K}} \psi_k u_k \Big\|_{V_\theta^B(\Gamma)}^p \leq 
C(N_0,p) \sum_{k \in \mathcal{K}} \Vert \psi_k u_k  \Vert_{V_\theta^B(\Gamma)}^p 
\end{equation}
for $u_k \in V_\theta^B(\Gamma)$.
\end{corollary}
\begin{lemma} \label{suml4}
We have that 
\begin{eqnarray} 
\label{sum5} \sum_{k \in \mathcal{K}} \Vert \psi_k g_k \Vert^q_{L_q(\Gamma)} & \leq & C^q \, \sum_{k \in \mathcal{K}} \Vert g_k \Vert^q_{L_q(\Gamma)}, \\  
\label{sum6} \sum_{k \in \mathcal{K}} \Vert \psi_k f \Vert^q_{L_q(\Gamma)} & \leq & (C N_0)^q \, \Vert f \Vert^q_{L_q(\Gamma)},
\end{eqnarray}
\begin{eqnarray} \label{sum9}
\sum_{k \in \mathcal{K}} [ \psi_k g_k ]^q_{q,r,\Gamma} & \leq & C^q \sum_{k \in \mathcal{K}} \Vert g_k \Vert^q_{W^r_q(\Gamma)} + \tilde{C} \, \lambda^{-q} \, \sum_{k \in \mathcal{K}} \Vert g_k \Vert^q_{L_\infty(\Gamma)}
\end{eqnarray}
and 
\begin{eqnarray} \label{sum10}
\sum_{k \in \mathcal{K}} [ \psi_k f ]^q_{q,r,\Gamma} & \leq & (C N_0)^q  \, \Vert f \Vert^q_{W^r_q(\Gamma)} + \tilde{C} \, \lambda^{-q-n+1} \, \Vert f \Vert^q_{L_\infty(\Gamma)}. 
\end{eqnarray}
\end{lemma}
{\bf Proof:}
Inequality (\ref{sum5}) is obvious. For (\ref{sum6}) note that
\[
\sum_{k \in \mathcal{K}} \sup_{x \in \Gamma} | \psi_k(x) | \leq C N_0.
\]
This implies 
\begin{eqnarray} \label{sum7}
& & \sum_{k \in \mathcal{K}} \int_{\Gamma} \int_{\Gamma} \frac{|\psi_k(x)|^q |f (x) - f (y)|^q}{|x-y|^{n-1+rq}} \; d\sigma(x)\;d\sigma(y) \nonumber \\
%& = & \int_{\Gamma} \sum_{l=1}^{N_0} \sum_{k \in \mathcal{K}} \int_{\Gamma \cap \Omega_k^l} \frac{|\psi_k(x)|^p |f (x) - f (y)|^p}{|x-y|^{n-1+rp}} \; d\sigma(x)\;d\sigma(y) \\
& \leq & (C N_0)^q \, \int_{\Gamma} \int_{\Gamma} \frac{|f (x) - f (y)|^q}{|x-y|^{n-1+rq}} \; d\sigma(x)\;d\sigma(y). 
\end{eqnarray}
Inequality (\ref{sum9}) follows from     
\begin{eqnarray} \label{sum8} 
[\psi_k]^q_{q,r,\Gamma} & = & \int_{\Gamma} \int_{\Gamma} \frac{|\psi_k (x)-\psi_k (y)|^q}{|x-y|^{n-1+rq}} \; d\sigma(x)\;d\sigma(y) \nonumber \\
& \leq & C^q \, \lambda^{-q} \, \int_{\mathbb{T}^{n-1}} \int_{\mathbb{T}^{n-1}} |x-y|^{-n+1+q(1-r)} \; d\sigma(x)\;d\sigma(y) \nonumber \\
& \leq & \tilde{C} \, \lambda^{-q}, 
\end{eqnarray}
and (\ref{sum10}) is obtained by combining (\ref{sum7}), (\ref{sum8}) and the fact that $|\mathcal{K}| \sim \lambda^{-n+1}$. \qed  
\begin{remark} \label{suml5}
{\rm In the same way as above one obtains    
\begin{eqnarray} \label{sum11}
\sum_{k \in \mathcal{K}} [ \partial_i (\psi_k g_k) ]^q_{q,r,\Gamma} & \leq & C^q \lambda^{-q} \sum_{k \in \mathcal{K}} \Vert g_k \Vert^q_{W^r_q(\Gamma)} + \tilde{C} \, \lambda^{-q-1} \, \sum_{k \in \mathcal{K}} \Vert \partial_i(g_k) \Vert^q_{L_\infty(\Gamma)} \nonumber \\
& & + \; \tilde{C} \, \lambda^{-q} \, \sum_{k \in \mathcal{K}} \Vert \partial_i g_k \Vert^q_{L_\infty(\Gamma)} + C^q \sum_{k \in \mathcal{K}} \Vert \partial_i g_k \Vert^q_{W^r_q(\Gamma)},
\end{eqnarray}
$i=1,\dots,n-1$. From this one concludes  
\begin{equation} \label{sum12}
\sum_{k \in \mathcal{K}} \Vert \psi_k u_k \Vert_{X_\theta^B(\Gamma)}^p \leq C(\lambda,p)T^\delta \sum_{k \in \mathcal{K}} \Vert u_k \Vert_{X_\theta^B(\Gamma)}^p + C^p \sum_{k \in \mathcal{K}} \Vert u_k \Vert_{X_\theta^B(\Gamma)}^p
\end{equation}
for $u_k \in \Mathring{X}_\theta^B(\Gamma)$, $p$ as in Section 4.1, and some $\delta > 0$.}
\end{remark}

\medskip

\paragraph*{Acknowledgements:} The research leading to this paper was carried 
out in part while the second author enjoyed the hospitality of the Institute of 
Applied Mathematics of Leibniz University Hannover. Moreover, we express our 
gratitude to E.V.\ Frolova, J. Seiler, and M.\ Wilke for helpful comments and 
discussions.

\end{document}